\documentclass[11pt]{article}
\usepackage{color,latexsym,amsmath,amssymb,path,amsthm,url}

\newcommand{\ignore}[1]{}

\newtheorem{theorem}{Theorem}[section]
\newtheorem{lemma}[theorem]{Lemma}
\newtheorem{corollary}[theorem]{Corollary}

\theoremstyle{definition}
\newtheorem{definition}[theorem]{Definition}
\newcommand{\Proof}[1]
        {
        \noindent
        \emph{Proof #1.}~
        }
\newsavebox{\smallProofsym}                     
\savebox{\smallProofsym}                        %
        {
        \begin{picture}(7,7)                    %
        \put(0,0){\framebox(6,6){}}             %
        \put(0,2){\framebox(4,4){}}             %
        \end{picture}                           %
        }                                       %

\newenvironment{theProof}[1]
        {
        \Proof{#1}}{\hfill\(\Box\)
        \medskip

        }
\usepackage{graphicx,psfrag}
\usepackage[rflt]{floatflt}

\newcommand{\placefig}[2]
        {\includegraphics[width=#2]{#1.eps}}
\usepackage{dsfont}

\newcommand{\RR}{\ensuremath{\mathbb R}}
\newcommand{\CC}{\ensuremath{\mathbb C}}

\newcommand{\fl}[1]{\mathsf{FL}_{#1}}

\def\nn{{\bf n}}

\def\deg{{\rm deg}}

\def\eps{{\varepsilon}}

\newcommand{\pts}{\mathcal P}
\newcommand{\circles}{\mathcal C}
\newcommand{\plane}{\Pi}

\begin{document}
\pagenumbering{arabic}

\title{Improved bounds for incidences between points and circles\thanks{%
Work on this paper by the first two authors was partially
supported by Grant 338/09 from the Israel Science Fund
and by the Israeli Centers of Research Excellence (I-CORE)
program (Center No.~4/11). Work by the first author was also
supported by NSF Grant CCF-08-30272, by Grant 2006/194
from the U.S.-Israel Binational Science Foundation, and
by the Hermann Minkowski-MINERVA Center for Geometry
at Tel Aviv University. Work on this paper by the third author
was partially supported by the Department of Defense through
a National Defense Science \& Engineering Graduate (NDSEG) Fellowship.
Work on this paper was also partially supported by the
ERC Advanced Grant No.~267165.
A preliminary version of this paper has appeared in
\emph{Proc. 29th ACM Annu. Symposium on Computational Geometry}, 2013.}}


\author{
Micha Sharir\thanks{%
School of Computer Science, Tel Aviv University,
Tel Aviv 69978, Israel.
{\sl michas@tau.ac.il} }
\and
Adam Sheffer\thanks{%
School of Computer Science, Tel Aviv University,
Tel Aviv 69978, Israel.
{\sl sheffera@tau.ac.il} }
\and
Joshua Zahl\thanks{%
Department of Mathematics, UCLA, Los Angeles, CA, 90024.
{\sl jzahl@mit.edu}
}}

\maketitle

\begin{abstract}
We establish an improved upper bound for the number of incidences
between $m$ points and $n$ circles in three dimensions.
The previous best known bound, originally established for the planar
case and later extended to any dimension $\ge 2$, is
$O^*\left( m^{2/3}n^{2/3} + m^{6/11}n^{9/11}+m+n \right)$, where the
$O^*(\cdot)$ notation hides polylogarithmic factors. Since all the points
and circles may lie on a common plane (or sphere), it is impossible to
improve the bound in $\RR^3$ without first improving it in the plane.

Nevertheless, we show that if the set of circles is required to
be ``truly three-dimensional'' in the sense that no sphere or plane
contains more than $q$ of the circles, for some $q\ll n$,
then for any $\eps>0$ the bound can be improved to
\[O\big(m^{3/7+\eps}n^{6/7} + m^{2/3+\eps}n^{1/2}q^{1/6} +
m^{6/11+\eps}n^{15/22}q^{3/22} + m + n\big). \] \vspace{-5.5mm}

\noindent For various ranges of parameters (e.g., when $m=\Theta(n)$
and $q = o(n^{7/9})$), this bound is smaller than the lower bound $\Omega^*(m^{2/3}n^{2/3}+m+n)$, which holds in two dimensions.

We present several extensions and applications of the new bound:
(i) For the special case where all the circles have the same radius, we
obtain the improved bound
$O\left(m^{5/11+\eps}n^{9/11} + m^{2/3+\eps}n^{1/2}q^{1/6} + m + n\right)$.
(ii) We present an improved analysis that removes the subpolynomial
factors from the bound when $m=O(n^{3/2-\eps})$ for any fixed $\varepsilon >0$.
(iii) We use our results to obtain the improved bound $O(m^{15/7})$ for the
number of mutually similar triangles determined by any set of $m$ points in $\RR^3$.

Our result is obtained by applying the polynomial partitioning technique
of Guth and Katz using a constant-degree partitioning polynomial (as was
also recently used by Solymosi and Tao).  We also rely on various
additional tools from analytic, algebraic, and combinatorial geometry.
\end{abstract}

\section{Introduction}

Recently, Guth and Katz \cite{GK11} presented the \emph{polynomial partitioning
technique} as a major technical tool in their solution of the famous planar
distinct distances problem of Erd\H os \cite{erd46}.  This problem can be
reduced to an incidence problem involving points and lines in $\RR^3$
(following the reduction that was proposed in \cite{ES11}), which can be
solved by applying the aforementioned polynomial partitioning technique.
The Guth--Katz result prompted various other incidence-related studies that
rely on polynomial partitioning (e.g., see \cite{BasuSombra14,FPSSZ14,Guth14,KMSS11,KMS11,SSo14,ST11,Zahl11}).
One consequence of these studies is that they have led to further developments
and enhancements of  the technique itself (evidenced  for example by the use of
induction in \cite{ST11}, and the use of two partitioning polynomials in
\cite{KMSS11,Zahl11}).  Also, the technique was recently applied to some
problems that are not incidence related: It was used to provide an alternative proof of the existence
of spanning trees with small crossing number in any dimension \cite{KMS11},
and to obtain improved algorithms for range searching with semialgebraic sets
\cite{AMS12}. Thus, it seems fair to say that applications and enhancements
of the polynomial partitioning technique form an active contemporary area of research in
combinatorial and computational geometry.

In this paper we study incidences between points and circles in three dimensions.
Let $\pts$ be a set of $m$ points and $\circles$ a set of $n$ circles in $\RR^3$.
We denote the number of point-circle incidences in $\pts\times \circles$ as
$I(\pts,\circles)$.  When the circles have arbitrary radii, the current best
bound for any dimension $d \ge 2$ (originally established for the planar case in \cite{ANPPSS04,AS02,MT06}, and later extended to higher dimensions by Aronov, Koltun, and Sharir \cite{AKS05})
is
\begin{equation}\label{eq:oldBound}
I(\pts,\circles) = O^*\left( m^{2/3}n^{2/3} + m^{6/11}n^{9/11}+m+n \right).
\end{equation}
The precise best known upper bound is
$O(m^{2/3}n^{2/3} + m^{6/11}n^{9/11}\log^{2/11}(m^3/n)+m+n)$ \cite{MT06}.

Since the three-dimensional case also allows $\pts$ and $\circles$ to
lie on a single common plane or sphere\footnote{
  There is no real difference between the cases of coplanarity and
  cosphericality of the points and circles, since the latter case can be
  reduced to the former (and vice versa) by means of the stereographic
  projection.},
the point-circle incidence bound in $\RR^3$ cannot be improved without first
improving the planar bound \eqref{eq:oldBound} (which is still an open problem, for about 10 years).
Nevertheless, as we show in this paper, an improved bound can be
obtained if the configuration of points and circles is ``truly
three-dimensional''
in the sense that no sphere or plane contains too many circles from
$\circles$. (Guth and Katz \cite{GK11} impose a similar assumption on the
maximum number of lines that can lie in a common plane or regulus; see also Sharir and Solomon~\cite{SSo14}.)
Our main result is given in the following theorem.
\begin{theorem} \label{th:arb}
Let $\pts$ be a set of $m$ points and let $\circles$ be a set of $n$
circles in $\RR^3$, let $\eps$ be an arbitrarily small positive constant,
and let $q\leq n$ be an integer.
If no sphere or plane contains more than $q$ circles of $\circles$, then
\begin{equation*}
I(\pts,\circles) = O\left(m^{3/7+\eps}n^{6/7} +
m^{2/3+\eps}n^{1/2}q^{1/6} + m^{6/11+\eps}n^{15/22}q^{3/22}
+ m + n\right),
\end{equation*}
where the constant of proportionality depends on $\eps$.
\end{theorem}
\paragraph{Remarks.} \hspace{-3mm} (1) In the planar case, the best
known lower bound for the number of point-circle incidences is
$\Omega^*(m^{2/3}n^{2/3}+m+n)$ (e.g., see \cite{PS04}).
Theorem \ref{th:arb} implies that for certain ranges of $m,n,$ and
$q$, a smaller upper bound holds in $\RR^3$.
This is the case, for example, when $m=\Theta(n)$ and $q = o(n^{7/9})$.
\medskip

\noindent (2) When $m>n^{3/2}$, we have $m^{3/7}n^{6/7}<m$ and
$m^{6/11}n^{15/22}q^{3/22}<m^{2/3}n^{1/2}q^{1/6}$.
Hence, we have
\[I(\pts,\circles)=O(m^{2/3+\eps}n^{1/2}q^{1/6}+m^{1+\eps}).\]
In fact, as the analysis in this paper will show, the first term in this bound
only arises from bounding incidences on certain potentially ``rich''
planes or spheres.
For $q=O(m^2/n^3)$ we have $I(\pts,\circles)=O(m^{1+\eps})$.
Informally, for such small values of $q$, when the number of points
is $\Omega(n^{3/2})$, a typical point can lie in only a small
(nearly constant) number of circles.
\medskip

\noindent (3) When $m\le n^{3/2}$, any of the terms except for $m$ can
dominate the bound.  However, if in addition
$q=O\left(\left(\frac{n^3}{m^2}\right)^{3/7}\right)$ then the bound becomes
\[ I(\pts,\circles)=O(m^{3/7+\eps}n^{6/7}+n^{1+\eps}),\]
for any $\varepsilon>0$.
Note also that the interesting range of parameters is
$m=\tilde\Omega(n^{1/3})$ and $m=\tilde O(n^2)$, where the $\tilde\Omega(\cdot)$ and $\tilde O(\cdot)$ notations hide terms of the form $m^\eps$ and $n^{\eps}$;
in the complementary ranges both the old and new bounds become (almost) linear
in $m+n$.  In the interesting range, the new bound is asymptotically
smaller than the planar bound given in \eqref{eq:oldBound} for $q$
sufficiently small (e.g., when $q=O\left( \left(\frac{n^3}{m^2}\right)^{3/7}\right)$ as above), and as noted, it is also smaller than the best known worst-case lower bound in the planar case for certain ranges of $m$ and $n$.

\noindent (4) Interestingly, the ``threshold'' value $m=\Theta(n^{3/2})$ where a quantitative change in the bound takes place (as noted in Remarks (2) and (3) above) also arises
in the study of incidences between points and lines in $\RR^3$
\cite{EKS11,GK10,GK11}.  See Section~\ref{furtherApplicationsSection}
for  a discussion of this threshold phenomenon.
\medskip

The proof of Theorem \ref{th:arb} is based on the polynomial
partitioning technique of Guth and Katz~\cite{GK11}, where we use a
constant-degree partitioning polynomial
in a manner similar to that used by Solymosi and Tao~\cite{ST11}.
(The use of constant-degree polynomials and the inductive arguments it
leads to are essentially the only similarities with the technique
of \cite{ST11}, which does not apply to circles in any dimension
since it cannot handle situations where arbitrarily many curves can
pass between any specific pair of points. Constant degree partitioning polynomials were also recently used in \cite{Guth14,SSo14}.)
The application of this
technique to incidences involving circles leads to new problems,
involving the handling of points that are incident to many circles
that are entirely contained in the zero set of the partitioning polynomial.
To handle this situation we turn these circles into lines using an
inversion transformation.  We then analyze the geometric and algebraic
structure of the transformed zero set using a variety of tools such
as \emph{flecnode polynomials} (as used in  \cite{GK11}), additional
classical $19^{\textrm{th}}$-century results in analytic geometry from \cite{sal15} (mostly related to \emph{ruled surfaces}), a very recent technique for analyzing surfaces that are ``ruled" by lines and circles \cite{NS12},
and some traditional tools from combinatorial geometry.
We note that while our results refer to the real affine case, part of our analysis
considers the setup in complex projective spaces, in which more classical results
from algebraic geometry can be brought to bear. By transporting the results back
to the real affine case, we obtain the properties that we wish to establish. See Section \ref{sec:reg} for details.
\vspace{2mm}

\noindent {\bf Removing the epsilons.} One disadvantage of the current use of constant-degree partitioning polynomials is that it leads to an upper bound
involving $\varepsilon$'s in some of the exponents (as stated in Theorem~\ref{th:arb}).
In Section  \ref{ssec:eps} we present an alternative approach, which uses partitioning polynomials of higher degree but requires a more involved analysis, for partially removing these $\varepsilon$'s. It yields the following theorem:
\begin{theorem} \label{th:noEps}
Let $\pts$ be a set of $m$ points and let $\circles$ be a set of $n$
circles in $\RR^3$, let $q\leq n$ be an integer, and let $m=O(n^{3/2-\delta})$, for some fixed arbitrarily small
constant $\delta>0$.
If no sphere or plane contains more than $q$ circles of $\circles$, then
\begin{equation} \label{mainBound}
I(\pts,\circles) \le A_{m,n}\left(m^{3/7}n^{6/7} +
m^{2/3}n^{1/2}q^{1/6} + m^{6/11}n^{15/22}q^{3/22}\log^{2/11}m + m + n\right),
\end{equation}
where $A_{m,n} = A^{\left\lceil \frac32 \cdot \frac{\log(m/n^{1/3})}
{\log(n^{3/2}/m)} \right\rceil +1}$, for some absolute constant $A>1$.
\end{theorem}
%
Note that $A_{m,n}$ grows slowly with the quantity $\log(m/n)$. For example, it is
$A$ for $m\le n^{1/3}$, $A^2$ for $n^{1/3}<m\le n^{4/5}$, and $A^3$ for
$n^{4/5}<m\le n$.

Recently, several other situations in which such $\varepsilon$'s can be removed were described in \cite{Zahl12}.
Our alternative technique seems to be sufficiently general, and we hope that appropriate variants of it could yield similar improvements of other bounds
that were obtained with constant-degree partitioning polynomials, such as the ones in \cite{ST11}.
\vspace{2mm}

\noindent {\bf Unit circles.}
In the special case where all the circles of $\circles$ have the
same radius, we derive the following improved bound.
\begin{theorem} \label{th:uni}
Let $\pts$ be a set of $m$ points and let $\circles$ be a set of $n$
unit circles in $\RR^3$, let $\eps$ be an arbitrarily small positive constant,
and let $q\leq n$ be an integer.  If no plane or sphere contains more than $q$
circles of $\circles$, then
\[
I(\pts,\circles) = O\left(m^{5/11+\eps}n^{9/11} +
m^{2/3+\eps}n^{1/2}q^{1/6} + m + n\right),
\]
where the constant of proportionality depends on $\eps$.
\end{theorem}

This improvement is obtained through the following steps. (i) We use the improved planar (or spherical) bound
$O(m^{2/3}n^{2/3}+m+n)$ for incidences with coplanar (or co-spherical) unit circles
(e.g., see~\cite{Sz}). (ii) We show that the number
of unit circles incident to at least three points in a given
set of $m$ points in $\RR^3$ is only $O(m^{5/2})$. (iii) We use this bound as a bootstrapping tool for deriving the
bound asserted in the theorem. The full details are presented in
Section \ref{sec:uni}.  Here too we can refine the bound of
Theorem \ref{th:uni} and get rid of the $\varepsilon$'s in the exponents,
for $m=O(n^{3/2-\varepsilon})$ for any $\varepsilon>0$.  The resulting refinement,
analogous to that of Theorem \ref{th:arb}, is given in Section \ref{sec:uni}.
\vspace{2mm}

\noindent {\bf An application: similar triangles.} Given a finite point
set $\pts$ in $\RR^3$ and a triangle $\Delta$, we denote by $F(\pts,\Delta)$
the number of triangles that are spanned by points of $\pts$ and are
similar to $\Delta$. Let $F(m)=\max_{|\pts|=m,\Delta}F(\pts,\Delta)$.
The problem of obtaining good bounds for $F(m)$ is motivated by
questions in exact pattern matching, and has been studied in several
previous works (see \cite{aaps07,ATT98,AS11,Br02}).
Theorem \ref{th:noEps} implies the bound
$F(m)=O(m^{15/7})$, which slightly improves upon the previous bound
of $O^*(m^{58/27})$ from \cite{AS11}; see also \cite{aaps07}.
The new bound is an almost immediate corollary of Theorem \ref{th:noEps},
while the previous bound requires a more complicated analysis.  This
application is presented in Section \ref{furtherApplicationsSection}.

\section{Algebraic preliminaries} \label{sec:pre}
We briefly review in this section the machinery needed for our analysis,
including the polynomial partitioning technique of Guth and Katz and
several basic tools from algebraic geometry.

\paragraph{Polynomial partitioning.}
In what follows, we regard the dimension $d$ of the ambient space as a (small) constant,
and we ignore the dependence on $d$ of the various constants of proportionality in the bounds.
Consider a set $\pts$ of $m$ points in $\RR^d$. Given a polynomial
$f \in \RR[x_1,\ldots, x_d]$, we define the \emph{zero set} of $f$ to be
$Z(f) = \{ p \in \RR^d \mid f(p)=0 \}$.  For $1 < r \le m$, we say that
$f \in \RR[x_1,\ldots, x_d]$ is an \emph{$r$-partitioning polynomial} for
$\pts$ if no connected component of $\RR^d \setminus Z(f)$ contains more than
$m/r$ points of $\pts$.  Notice that there is no restriction on the number of
points of $\pts$ that lie in $Z(f)$.

The following result is due to Guth and Katz \cite{GK11}.
A detailed proof can also be found in \cite{KMS11}.
\begin{theorem}\label{th:partition}{\bf (Polynomial partitioning \cite{GK11})}
Let $\pts$ be a set of $m$ points in $\RR^d$. Then for every $1< r \le m$,
there exists an $r$-partitioning polynomial $f \in \RR[x_1,\ldots, x_d]$ of
degree $O(r^{1/d})$. \hfill\(\Box\)
\end{theorem}
To use such a partitioning effectively, we also need a bound on the maximum
possible number of cells of the partition.  Such a bound is provided by the following theorem.
\begin{theorem}\label{th:Warren}{\bf (Warren's theorem \cite{Warren68}; cf.\ also \cite{Alon95})}
For a polynomial $f \in \RR[x_1,\ldots, x_d]$ of degree $k$,
the number of connected components of $\RR^d \setminus Z(f)$ is $\displaystyle O\left((2k)^d\right)$. \hfill\(\Box\)
\end{theorem}
Consider an $r$-partitioning polynomial $f$ for a point-set $\pts$, as
provided in Theorem \ref{th:partition}.  The number of cells in the partition
is equal to the number of connected components of $\RR^d \setminus Z(f)$, which, by Theorem \ref{th:Warren}, is $O((r^{1/d})^d) = O(r)$ (recall that
$f$ is of degree $O(r^{1/d})$ and that $d$ is treated as a fixed constant --- 3 in our case).  It follows that the bound on the number of
points in each cell, namely $m/r$, is asymptotically best possible.

We will also rely on the following classical result, somewhat similar to Warren's theorem (the following formulation is taken from \cite{Alon95}).
\begin{theorem}\label{th:Milnor}{\bf (Milnor--Thom's theorem \cite{Mil64,Thom65})}
Let $V$ be a real variety in $\RR^d$, which is the solution set of the real polynomial equations
\[ f_i(x_1,\ldots,x_d) = 0 \qquad (i = 1,\ldots,m), \]
and suppose that the degree of each polynomial $f_i$ is at most $k$. Then the number of connected
components of $V$ is at most $k(2k-1)^{d-1}$. \hfill\(\Box\)
\end{theorem}

Since this paper studies incidences in a three-dimensional space, we will
only apply the above theorems for $d=3$.

\paragraph{Real and complex varieties.}
Let $F=\RR$ or $\CC$ be the underlying field, and let $I\subset F[x_1,x_2,x_3]$ be an ideal. We define $Z(I)\subset F^3$ to be the variety corresponding to $I$; this is the common vanishing locus of all the polynomials $f\in I$. When it is not clear from the context whether $F=\RR$ or $\CC$, we resolve the ambiguity by writing $Z_{\RR}(I)$ or $Z_{\CC}(I)$, respectively. If $f\in F[x_1,x_2,x_3]$, we will abuse notation and write $Z(f)$ instead of $Z((f))$ (where $(f)$ is the ideal generated by $f$). If $Z\subset F^3$ is a variety, we define $\mathbf{I}(Z)\subset F[x_1,x_2,x_3]$ to be the ideal of polynomials vanishing on $Z$ (we will use bold typeface $\mathbf{I}$ to distinguish it from the inversion transform $I\colon\RR^3\to\RR^3$ that will appear in Section \ref{sec:reg}).

Let $R=\RR[x_1,\ldots,x_d]$, and let $I\subset R$ be an affine ideal. We say that $I$ is a \emph{real ideal} if for every set $g_1,\ldots,g_t\in R$, we have
\begin{equation*}
g_1^2+\ldots+g_t^2\in I\ \Longrightarrow g_i \in I\ \textrm{for}\ i=1,\ldots,t.
\end{equation*}

Intuitively, algebraic geometry over the reals tends to be more pathological than algebraic geometry over $\CC$. Many of the problematic cases do not occur if we only work with real ideals. For example, the following theorem presents two properties of real ideals.
\begin{theorem} \label{th:real}
(i) Let $J\subset \RR[x_1,\ldots,x_d]$ be an ideal. Then $J=\mathbf{I}(Z_{\RR}(J))$ if and only if $J$ is a real ideal.\\
(ii) Let $f\in \RR[x_1,\ldots,x_d]$ be an irreducible polynomial. Then $(f)$ is a real ideal if and only if $\dim Z(f) = d-1$.
\end{theorem}
A nice introduction to real ideals can be found in \cite{BCR98}. The two parts of the above theorem
are Theorem 4.1.4 and Theorem 4.5.1 of \cite{BCR98}, respectively.

Given a variety $Z \subset \RR^3$, the \emph{complexification} $Z^*\subset \CC^3$
of $Z$ is the smallest complex variety that contains $Z$ (in the
sense that any other complex variety that contains $Z$ also contains $Z^*$,
e.g., see \cite{RV02,Whit57}). As shown in \cite[Lemma 6]{Whit57},
such a complexification always exists, and $Z$ is precisely the locus of real
points of $Z^*$. More specifically, $Z^* = Z_\CC(\mathbf{I}(Z))$.

According to \cite[Lemma 7]{Whit57}, there is a bijection between the set of
irreducible components of $Z$ and the set of irreducible components of $Z^*$,
such that each real component is the real part of its corresponding complex component.
Specifically, the complexification of an irreducible variety is irreducible.

\paragraph{B\'ezout's theorem.}
We also need the following basic property of zero sets of polynomials in
the plane (for further discussion see \cite{CLO1,CLO2}).
\begin{theorem}\label{th:bezout}{\bf (B\'ezout's theorem)}
Let $f,g$ be two polynomials in $\RR[x_1,x_2]$ or $\CC[x_1,x_2]$ of degrees $D_f$ and $D_g$,
respectively.
(i) If $Z(f)$ and $Z(g)$ have a finite number of common points, then this
number is at most $D_f D_g$.
(ii) If $Z(f)$ and $Z(g)$ have an infinite number of (or
just more than $D_f D_g$) common points, then $f$ and $g$ have a common (nontrivial) factor. \hfill\(\Box\)
\end{theorem}
The following result, also used in the sequel, is somewhat related to B\'ezout's theorem,
and holds in complex projective spaces of any dimension (e.g., see \cite{Fu98}; for a formal
definition of three-dimensional complex projective space $\CC\mathbf{P}^3$, in which we will
apply the following theorem and other tools, see Section \ref{sec:reg}, and recall the comments made earlier concerning the passage between the real and complex setups).
\begin{theorem}\label{th:bezout3}
Let $Z_1$ and $Z_2$ be pure-dimensional varieties (every irreducible component of a pure-dimensional variety has the same dimension) in $d$-dimensional complex projective space, with $\operatorname{codim}Z_1 + \operatorname{codim}Z_2=d$. Then if $Z_1\cap Z_2$ is a zero-dimensional set of points, this set is finite. \hfill\(\Box\)
\end{theorem}

The following lemma is a consequence of Theorem~\ref{th:bezout}. Its proof
is given in Guth an Katz \cite[Corollary 2.5]{GK10} and in Elekes, Kaplan, and Sharir \cite[Proposition 1]{EKS11}.
\begin{lemma}\label{le:commonLines} {\bf (Guth and Katz \cite{GK10})}
Let $f$ and $g$ be two polynomials in $\RR[x_1,x_2,x_3]$ (or in $\CC[x_1,x_2,x_3]$) of respective degrees $D_f$
and $D_g$, such that $f$ and $g$ have no common factor. Then there are at
most $D_fD_g$ lines on which both $f$ and $g$ vanish identically. \hfill\(\Box\)
\end{lemma}

\paragraph{Flecnode polynomial.}
A \emph{flecnode} of a surface $Z$ in $\CC^3$ is a point $p\in Z$ for which there exists a line
that passes through $p$ and agrees with $Z$ at $p$ to order three.
That is, if $Z=Z(f)$
and the direction of the line is $v=(v_1,v_2,v_3)$ then
\begin{equation}\label{tangOrdThreeEqn}
f(p) = 0, \qquad \nabla_v f(p) = 0, \qquad  \nabla_v^2 f(p) = 0,
\qquad  \nabla_v^3 f(p) = 0,
\end{equation}
where $\nabla_v f, \nabla^2_v f, \nabla^3_v f$ are, respectively, the
first, second, and third-order derivatives of $f$ in the direction $v$.

If $Z=Z(f)$ is a surface in $\RR^3$, we say that $p\in Z$ is a flecnode of $Z$ if $p$ is a flecnode of the corresponding complex surface $(Z(f))^*$.


The \emph{flecnode polynomial} of $f$, denoted $\fl{f}$, is the polynomial
obtained by eliminating $v$ from the last three equations in \eqref{tangOrdThreeEqn}.  Note that the corresponding
polynomials of the system are homogeneous in $v$ (of respective degrees
$1$, $2$, and $3$).  We thus have a system of three equations in six variables.
Eliminating the variables $v_1,v_2,v_3$ results in a single polynomial equation in
$p=(x_1,x_2,x_3)$, which is the desired flecnode polynomial. By construction,
the flecnode polynomial of $f$ vanishes on all the flecnodes of $Z(f)$.
The following results, also mentioned in \cite[Section 3]{GK11}, are taken
from Salmon \cite[Chapter XVII, Section III]{sal15}.

\begin{lemma}\label{le:flec1}
Let $Z\subset\RR^3$ be a surface, with $Z=Z(f)$ for a polynomial $f \in \RR[x_1,x_2,x_3]$ of degree $d\ge 3$. Then $\fl{f}$ is a real polynomial (i.e., an element of $\RR[x_1,x_2,x_3]$), and it has degree at most $11d-24$. \hfill\(\Box\)
\end{lemma}

\begin{definition}
An algebraic surface $S$ in three-dimensional space (we restrict our attention to $\RR^3$,
$\CC^3$, and $\CC\mathbf{P}^3$) is said to be \emph{ruled} if
every point of $S$ is incident to a straight line that is fully
contained in $S$.  Equivalently, $S$ is a (two-dimensional) union of lines.\footnote{%
  We do not insist on the more restrictive definition used in
  differential (or in algebraic) geometry, which requires the ruling lines to form
  a smooth 1-parameter family; cf. \cite[Chapter III]{Beau96} and \cite[Section V.2]{Hart83}.}
We say that an irreducible surface $S$ is \emph{triply ruled} if for
every point on $S$ there are (at least) three straight lines
contained in $S$ and passing through that point.
As is well known,
the only triply ruled surfaces are planes (e.g., see \cite[Lecture 16]{FT07}; while this reference only provides proofs for the case of $\RR^3$, proofs for the cases of $\CC^3$ and $\CC\mathbf{P}^3$ are also known).
We say that an irreducible surface $S$ is \emph{doubly ruled}
if it is not triply ruled and for every point on $S$ there are
(at least) two straight lines contained in $S$ and passing
through that point. It is well known that the
only doubly ruled surfaces are the hyperbolic paraboloid and
the hyperboloid of one sheet (again, see \cite[Lecture 16]{FT07}).
Finally, we say that an irreducible ruled surface is
\emph{singly ruled} if it is neither doubly nor triply ruled.
\end{definition}

\begin{lemma}\label{le:flec3}
Let $Z\subset\RR^3$ be a surface with $Z=Z(f)$ for a polynomial $f \in \RR[x_1,x_2,x_3]$ of degree $d\ge 3$. Then every
line that is fully contained in $Z$ is also fully contained in $Z(\fl{f})$.
\end{lemma}
\begin{theProof}{\!\!}
Every point on any such line is a flecnode of $Z$, so $\fl{f}$ vanishes
identically on the line.
\end{theProof}
\begin{theorem}\label{th:flec2} {\bf (Cayley--Salmon \cite{sal15})}
Let $Z\subset\RR^3$ be a surface with $Z=Z(f)$ for a polynomial $f \in \RR[x_1,x_2,x_3]$ of degree $d\ge 3$.
Then $Z$ is ruled if and only if  $Z \subseteq Z(\fl{f})$. \hfill\(\Box\)
\end{theorem}

\begin{corollary} \label{co:ruled}
Let $Z\subset\RR^3$ be a surface with $Z=Z(f)$ for an irreducible polynomial $f \in \RR[x_1,x_2,x_3]$ of degree $d\ge 3$.
If $Z$ contains more than $d(11d-24)$ lines then $Z$ is a ruled surface.
\end{corollary}
\begin{theProof}{\!\!}
Lemma \ref{le:commonLines} and Lemma \ref{le:flec3} imply that in this case
$f$ and $\fl{f}$ have a common factor.  Since $f$ is irreducible, $f$ divides
$\fl{f}$, and Theorem \ref{th:flec2} completes the proof.
\end{theProof}

A modern treatment of the Cayley--Salmon theorem can be found in a more
recent work by Landsberg \cite{Land99}.  (The results in \cite{Land99} are
considerably more general, but we state here only the special case related
to the Cayley--Salmon theorem.)

\begin{theorem}\label{th:landsberg} {\bf (Landsberg \cite{Land99})}
Let $Z$ be a surface in $\CC^3$, and let $Z=Z(f)$ for a polynomial $f$ of degree $d\ge 3$. Then $Z$ is ruled if and only if $Z \subseteq Z(\fl{f})$. \hfill\(\Box\)
\end{theorem}

\section{The main theorem} \label{sec3}
In this section we prove Theorem \ref{th:arb}, which we restate for the
convenience of the the reader.

\paragraph{Theorem \ref{th:arb}}
\emph{Let $\pts$ be a set of $m$ points and let $\circles$ be a set of $n$
circles in $\RR^3$, let $\eps$ be an arbitrarily small positive constant,
and let $q\leq n$ be an integer.  If no sphere or plane contains more than $q$
circles of $\circles$, then
\[
I(\pts,\circles) = O\left(m^{3/7+\eps}n^{6/7} + m^{2/3+\eps}n^{1/2}q^{1/6} +
m^{6/11+\eps}n^{15/22}q^{3/22} + m + n\right),
\]
where the constant of proportionality depends on $\eps$.
}\vspace{3mm}

\hspace{-9.3mm} \begin{theProof}{\!\!}
The proof proceeds by induction on $m+n$.
Specifically, we prove by induction that, for any fixed $\eps>0$, there exist constants $\alpha_1,\alpha_2$ such that
\[
I(\pts,\circles) \le \alpha_1\left(m^{3/7+\eps}n^{6/7} +
m^{2/3+\eps}n^{1/2}q^{1/6} +
m^{6/11+\eps}n^{15/22}q^{3/22} \right) + \alpha_2(m + n).
\]
Let $n_0$ be a constant (whose concrete choice will be made later).
The base case where $m+n<n_0$ can be dealt with by choosing $\alpha_1$ and $\alpha_2$ sufficiently large.

We start by recalling a well-known simple, albeit weaker bound.
The incidence graph $G\subseteq \pts\times\circles$ whose edges are
the incident pairs in $\pts\times\circles$, cannot contain $K_{3,2}$
as a subgraph, because two circles have at most two intersection
points.  By the K\H ov\'ari--S\'os--Tur\'an theorem (e.g., see
\cite[Section 4.5]{mat02}),
$I(\pts,\circles) = |G| = O\left(n^{2/3}m+n\right)$.
This immediately implies the theorem if $m=O\left(n^{1/3}\right)$
(the resulting bound is $O(n)$ in this case).
Thus we may assume that $n=O\left(m^{3}\right)$.

We next apply the polynomial partitioning technique.  Specifically, we
set $r$ as a sufficiently large constant (whose value depends on
$\eps$ and will be determined later), and apply the polynomial partitioning
theorem (Theorem \ref{th:partition}) to obtain an $r$-partitioning
polynomial $f$.  According to the theorem, $f$ is of degree
$D = O\left(r^{1/3}\right)$ and $Z(f)$ partitions $\RR^3$ into
maximal connected cells, each containing at most $m/r$ points of
$\pts$.  As already noted, Warren's theorem (Theorem \ref{th:Warren})
implies that the number of cells is  $O(r)$.

Let $\circles_0$ denote the subset of circles of $\circles$ that are fully
contained in $Z(f)$, and let $\circles^\prime = \circles \setminus \circles_0$.
Similarly, set $\pts_0 = \pts\cap Z(f)$ and $\pts^\prime = \pts \setminus \pts_0$.
Notice that
\begin{equation} \label{eq:part}
I(\pts,\circles) = I(\pts_0,\circles_0)
+ I(\pts_0,\circles^\prime) + I(\pts^\prime,\circles^\prime).
\end{equation}
The terms $I(\pts_0,\circles^\prime)$ and
$I(\pts^\prime,\circles^\prime)$ can be bounded using techniques
(detailed below) that are by now fairly standard.
On the other hand, bounding $I(\pts_0,\circles_0)$ is the main
technical challenge in this proof.
Other works that have applied the polynomial partitioning technique,
such as \cite{KMSS11,KMS11,ST11,Zahl11, Zahl12}, also spend most of
their efforts on incidences with
curves that are fully contained in the zero set of the partitioning
polynomial (where these curves are either original input curves specified in the statement of the problem, or the
intersections of
input surfaces with the zero set of a partitioning polynomial).

\paragraph{Bounding $I(\pts_0,\circles^\prime)$ and
$I(\pts^\prime,\circles^\prime)$.}
For a circle $C \in \circles^\prime$, let $\plane_C$ be the plane that
contains $C$, and let $f_C$ denote the restriction of $f$ to $\plane_C$.
Since $C$ is not contained in $Z(f_C)$, $f_C$ and the irreducible
quadratic equation of $C$ within $\plane_C$ do not have any common factor.
Thus by B\'ezout's theorem (Theorem \ref{th:bezout}), $C$ and $Z(f_C)$
have at most $2\cdot \deg(f_C) = O\left(r^{1/3}\right)$ common
points. This immediately implies
\vspace{-4mm}

\begin{equation} \label{eq:c'p0}
I(\pts_0,\circles^\prime) = O\left(nr^{1/3}\right).
\end{equation}

Next, let us denote the cells of the partition as $K_1, \ldots, K_s$
(recall that $s=O(r)$ and that the cells are open).  For $i=1,\ldots,s$, put
$\pts_i=\pts\cap K_i$ and let $\circles_i$ denote the set of circles
in $\circles^\prime$ that intersect $K_i$.  Put $m_i = |\pts_i|$ and
$n_i=|\circles_i|$, for $i=1,\ldots, s$. Note that $|\pts^\prime|=\sum_{i=1}^{s}m_i,$ and recall
that $m_i \le m/r$ for every $i$.  The above bound of
$O\left(r^{1/3}\right)$ on the number of intersection points
of a circle $C\in \circles^\prime$ and $Z(f)$ implies that each circle
enters $O\left(r^{1/3}\right)$ cells (a circle has to intersect
$Z(f)$ when moving from one cell to another).
This implies $\sum_{i=1}^s n_i = O\left(nr^{1/3}\right)$.

Notice that $I(\pts^\prime,\circles^\prime) = \sum_{i=1}^sI(\pts_i,\circles_i)$,
so we proceed to bound the number of incidences within a cell $K_i$.
From the induction hypothesis, we get

\begin{align}
  I(\pts^\prime,\circles^\prime) & \le
  \sum_{i=1}^{s}\left(\alpha_1\left(m_i^{3/7+\eps}n_i^{6/7} +
  m_i^{2/3+\eps}n_i^{1/2}q^{1/6} + m_i^{6/11+\eps}n_i^{15/22}q^{3/22}\right)
  + \alpha_2(m_i + n_i)\right) \nonumber\\
 & \le\sum_{i=1}^{s}
  \left(\alpha_1\left(\left(\frac{m}{r}\right)^{3/7+\eps}n_i^{6/7} +
  \left(\frac{m}{r}\right)^{2/3+\eps}n_i^{1/2}q^{1/6}
+ \left(\frac{m}{r}\right)^{6/11+\eps}n_i^{15/22}q^{3/22}\right)\right) \nonumber \\
& \qquad\qquad\qquad\qquad\qquad\qquad\qquad\qquad\qquad\quad
  +\alpha_2\left(|\pts'|+ \sum_{i=1}^{s}n_i \right). \label{eq:c'p'1}
\end{align}
Since $\sum_{i=1}^s n_i = O\left(nr^{1/3}\right)$, H\"older's inequality implies
\begin{equation}\label{eq:hol}
\begin{split}
\sum_{i=1}^{s}n_i^{6/7} &= O\left(\left(nr^{1/3}\right)^{6/7} \cdot
r^{1/7} \right) = O\left(n^{6/7}r^{3/7}\right),\\
\sum_{i=1}^{s}n_i^{1/2} &= O\left(\left(nr^{1/3}\right)^{1/2} \cdot
r^{1/2} \right) = O\left(n^{1/2}r^{2/3}\right), \\
\sum_{i=1}^{s}n_i^{15/22} &= O\left(\left(nr^{1/3}\right)^{15/22} \cdot
r^{7/22} \right) = O\left(n^{15/22}r^{6/11}\right).
\end{split}
\end{equation}
By combining \eqref{eq:c'p'1} and \eqref{eq:hol}, we obtain
\begin{equation*}
 \begin{split}
  I(\pts^\prime,\circles^\prime) \le
  \alpha_1\cdot \frac{c\left(m^{3/7+\eps}n^{6/7} +
  m^{2/3+\eps}n^{1/2}q^{1/6}+m^{6/11+\eps}n^{15/22}q^{3/22}\right)}{r^{\eps}}\\
  \qquad\qquad\qquad\qquad\qquad\qquad\qquad\qquad\qquad\qquad\quad+
  \alpha_2\left(|\pts'|+ cnr^{1/3}\right),
 \end{split}
\end{equation*}
for a suitable constant $c>0$. Notice that the bound in \eqref{eq:c'p0} is proportional to the last term in this bound, and that this term is
dominated by $O(m^{3/7}n^{6/7})$ since we assume that $n=O\left(m^{3}\right)$ and that $r$ is constant.
Choosing $r$ to be sufficiently large, so that $r^\varepsilon > 4c$, and choosing $\alpha_1\gg\alpha_2r^{1/3}$, we can ensure that
\begin{equation} \label{eq:c'p'2}
I(\pts_0 \cup \pts^\prime,\circles^\prime) \le
\frac{\alpha_1}{3}\left(m^{3/7+\eps}n^{6/7} +
m^{2/3+\eps}n^{1/2}q^{1/6}+m^{6/11+\eps}n^{15/22}q^{3/22}\right)
+\alpha_2 |\pts'|.
\end{equation}

\paragraph{Bounding $I(\pts_0,\circles_0)$: Handling shared points.}
We are left with the task of bounding the number of incidences between the set $\pts_0$ of points
of $\pts$ that are contained in $Z(f)$ and the set $\circles_0$ of circles of $\circles$ that are
fully contained in $Z(f)$.  We call a point of $\pts_0$ \emph{shared} if it
is contained in the zero sets of at least two distinct irreducible factors of $f$,
and otherwise we call it \emph{private}.  We first consider the case of shared points; in this argument it is simplest to work over $\CC$.

Let $\pts_s$ denote the subset of points in $\pts_0$ that are shared,
and put $m_{s}=|\pts_s|$. Let $f_1$ be the square-free (over $\CC$) part of $f,$ so in particular $Z_{\CC}(f)=Z_{\CC}(f_1)$. If $p\in\pts_s$ is a shared point, then $p$ lies in at least two distinct irreducible (over $\RR$) components of $Z_{\RR}(f)$, and thus $p$ lies in at least two irreducible (over $\CC$) components of $Z_{\CC}(f_1)$. Thus $p$ lies in the singular set of $Z_{\CC}(f_1),$ and in particular $p\in Z_{\CC}(f_1)\cap Z_{\CC}(f_1^\prime)$, where $f_1^\prime = e\cdot\nabla f_1$ and $e$ is a generically chosen unit vector (i.e., $f_1^\prime$ is a partial derivative of $f_1$ in a generic direction). Note that $\deg f_1, \deg f_1^\prime \leq D$, so $\gamma = Z_{\CC}(f_1)\cap Z_{\CC}(f_1^\prime)$ is an algebraic space curve of degree at most $D^2$ (e.g., see \cite[Exercise 11.6]{harris}). The curve $\gamma$ contains at most $D^2$ irreducible components, and thus $\gamma$ contains at most $D^2$ (complex) circles. We conclude that there are at most $D^2m_s$ incidences between points from $\pts_s$ and circles whose complexification is contained in $Z_{\CC}(f_1)\cap Z_{\CC}(f_1^\prime)$.

It remains to bound the number of incidences between points in
$\pts_s$ and circles of $\circles_0$ whose complexification is not contained in $Z_{\CC}(f_1)\cap
Z_{\CC}(f_1^\prime)$ (that is, circles that are contained in $Z_{\RR}(f)$ but whose complexification is not fully contained in
$Z_{\CC}(f_1^\prime)$). Let $C$ be a circle whose complexification $C^*$ is not contained in $Z_{\CC}(f_1^\prime)$, and let $\Pi\subset\CC^3$ be the 2--plane containing $C^*$. If we identify $\Pi$ with $\CC^2$, then the restriction of $f_1^\prime$ to $\Pi$ is a polynomial $\tilde f_1^\prime\in \CC[x_1,x_2]$. By B\'ezout's theorem (Theorem \ref{th:bezout}), $C^*$ and $Z_{\CC}(\tilde f_1^\prime)$ intersect in at most $2D$ points. Thus $C^*$ and $Z_{\CC}(f_1^\prime)$ intersect in at most $2D$ points. This in turn
implies that $|C\cap\pts_s|\le 2D$. Therefore, by taking $\alpha_2$ (and consequently also $\alpha_1$) to be sufficiently large, we have

\begin{equation}\label{eq:shared}
I(\pts_s,\circles_0) \le \frac12 {D^2}m_s + 2D n \le \alpha_2 (m_s + n/3).
\end{equation}

\paragraph{Bounding $I(\pts_0,\circles_0)$: Handling private points.}
Let $\pts_p = \pts_0 \setminus \pts_s$ denote the set of private
points in $\pts_0$.
Recall that each private point is contained in the zero set of a
single irreducible factor of $f$.
Let $f_1,f_2,\ldots,f_t$ be the factors of $f$ whose zero sets are
planes or spheres.
For $i=1,\ldots,t$, set $\pts_{p,i}^{(1)}=\pts_p \cap  Z(f_i)$  and
$m_{p,i}=|\pts_{p,i}^{(1)}|$.
Put $\pts_p^{(1)}=\bigcup_{i=1}^{t} \pts_{p,i}^{(1)}$ and $m_p^{(1)}= |\pts_p^{(1)}| = \sum_{i=1}^t m_{p,i}$.
Let $n_{p,i}$ denote the number of circles of $\circles_0$ that are
fully contained in $Z(f_i)$.
Notice that (i) $t \le D = O\left(r^{1/3}\right)$, (ii) $n_{p,i} \le q$
for every $i$, and (iii) $\sum_{i=1}^t n_{p,i} \le n$ (we may ignore circles
that are fully contained in more than one component, since these will not
have incidences with private points).  Applying \eqref{eq:oldBound} and using
the fact that there are no hidden polylogarithmic terms in
the linear part of \eqref{eq:oldBound}, we obtain\footnote{%
  Notice that the dependency of this bound on $n$ and $q$ is better than
  the one in the bound of the theorem. This latter worse bound is the one
  that is preserved under the partition-based induction.}
\begin{align*}
I(\pts_p^{(1)},\circles_0) &=
\sum_{i=1}^t \left(O^*\left( m_{p,i}^{2/3}n_{p,i}^{2/3} +
m_{p,i}^{6/11}n_{p,i}^{9/11}\right)+O(m_{p,i}+n_{p,i})\right)  \\
&= \sum_{i=1}^t \left( O^*\left(
m_{p,i}^{2/3}n_{p,i}^{1/3}q^{1/3} +
m_{p,i}^{6/11}n_{p,i}^{5/11}q^{4/11}\right)+ O(m_{p,i}+n_{p,i})\right)  \\
&= O^*\left( m^{2/3}n^{1/3}q^{1/3} +
m^{6/11}n^{5/11}q^{4/11}\right)+O\left(m_p^{(1)}+n\right),
\end{align*}
where the last step uses H\"older's inequality;
it bounds (twice) $\sum_i m_{p,i} = m_p^{(1)}$
by $m$.  Since $q\le n$, it follows that when $\alpha_1$
and $\alpha_2$ are sufficiently large, we have
\begin{equation}\label{eq:reg1}
I(\pts_p^{(1)},\circles_0) \le
\frac{\alpha_1}{3}\left( m^{2/3+\eps}n^{1/2}q^{1/6} +
m^{6/11+\eps}n^{15/22}q^{3/22} \right)+\alpha_2(m_p^{(1)}+n/3).
\end{equation}

Let $\pts_p^{(2)}=\pts_p\setminus \pts_p^{(1)}$ be the set of private
points that lie on the zero sets of factors of $f$ that are neither
planes nor spheres, and put $m_p^{(2)}=\big|\pts_p^{(2)}\big|$.
To handle incidences with these points we require the following lemma,
which constitutes a major component of our analysis and which is proved
in Section \ref{sec:reg} (somewhat similar results can be found in
\cite{Hwang05,Land03}). First, a definition.
\begin{definition}
Let $g$ be an irreducible polynomial in $\RR[x_1,x_2,x_3]$ such that
$Z(g)$ is a 2-dimensional surface.  We say that
a point $p \in Z(g)$ is \emph{popular} if it is incident to at least
$44(\operatorname{deg}g)^2$ circles that are fully contained in $Z(g)$.
\end{definition}
\begin{lemma} \label{le:reg}
An irreducible algebraic surface that is neither a plane nor a sphere
cannot contain more than two popular points.
\end{lemma}
The lemma implies that the number of incidences between popular points
of $\pts_p^{(2)}$ (within their respective irreducible components of
$Z(f)$, whose number is at most $D/2$) and circles of $\circles_0$ is at most
$2(D/2)n=Dn \le \alpha_2 n/3$ (the latter inequality holds if $\alpha_2$
is chosen sufficiently large with respect to $D$).  The number of incidences between non-popular
points of $\pts_p^{(2)}$ and circles of $\circles_0$ is at most
$m_p^{(2)}\cdot 44D^2 \le \alpha_2m_p^{(2)}$
(again for a sufficiently large value of $\alpha_2$).
Combining this with \eqref{eq:part}, \eqref{eq:c'p'2}, \eqref{eq:shared},
and \eqref{eq:reg1}, we get
\[
I(\pts,\circles) \le \alpha_1\left(m^{3/7+\eps}n^{6/7} + m^{2/3+\eps}n^{1/2}q^{1/6} +
m^{6/11+\eps}n^{15/22}q^{3/22} \right) + \alpha_2(m + n).\]
This establishes the induction step, and thus completes the proof of the theorem.
\end{theProof}

\noindent{\bf Remarks.}
(1) Note that we have actually shown that
$$
I(\pts_0,\circles_0) = O\left( mD^2+nD \right),
$$
regardless of the degree $D$ of $f$. It is the term $mD^2$ that becomes
too large when $D$ itself is too large. In this part of the analysis we
addressed this issue by taking $D$ to be a constant. In the refined analysis
given in the next subsection we use non-constant, albeit still small, values
for $D$, thereby slightly refining the bound.
\medskip

\noindent
(2) To see why $m^{3/7+\eps}n^{6/7}$ is the best choice for the leading term,
let us denote the leading term as $m^{a+\eps}n^b$ and observe the following
restrictions on $a$ and $b$:
(i) For $r$ to cancel itself in the analysis of the cells of the partition (up to a power
of $\eps$), we require $a\ge 1-2b/3$. (ii) For $n=O(m^3)$ to imply
$n=O(m^an^b)$, we must have $a+3b\ge3$. Combining both constraints,
with equalities, results in the term $m^{3/7+\eps}n^{6/7}$.
\medskip

\noindent (3)
We believe that the terms in our bounds that depend on $q$ can be significantly
improved, by using a more careful analysis. The exponents in these terms are chosen
so as to make them satisfy the induction step. However, in doing so,
in each cell of the partition,
we use the same bound $q$ on the maximum number of coplanar or cospherical circles
among those that cross the cell. Since the number of circles that cross a cell goes down, the bound $q$ should also
decrease. We do not know how to handle this issue rigorously, and leave it as an open problem for further research.

\subsection{Removing the epsilons} \label{ssec:eps}
In this section we will show that, for any $\delta >0$, when
$n=O(m^{3/2-\delta})$, the epsilons from the bound of Theorem \ref{th:arb}
can be removed. This is what Theorem \ref{th:noEps} asserts; we repeat its statement
for the convenience of the reader.

\noindent {\bf Theorem \ref{th:noEps}.} \emph{Let $\pts$ be a set of $m$ points and let $\circles$ be a set of $n$
circles in $\RR^3$, let $q\leq n$ be an integer, and let $m=O(n^{3/2-\delta})$, for some fixed arbitrarily small
constant $\delta>0$.
If no sphere or plane contains more than $q$ circles of $\circles$, then
\begin{equation*}
I(\pts,\circles) \le A_{m,n}\left(m^{3/7}n^{6/7} +
m^{2/3}n^{1/2}q^{1/6} + m^{6/11}n^{15/22}q^{3/22}\log^{2/11}m + m + n\right),
\end{equation*}
where $A_{m,n} = A^{\left\lceil \frac32 \cdot \frac{\log(m/n^{1/3})}
{\log(n^{3/2}/m)} \right\rceil+1}$, for some absolute constant $A>1$.}
\vspace{2mm}

\hspace{-9mm} \begin{theProof}{\!\!}
We define $\pts$, $\pts_0$, $\circles$, $\circles'$, etc., as in the proof of Theorem \ref{th:arb}.
The proof is similar to the one of Theorem \ref{th:arb}, except that it works in stages, so that in each stage we enlarge the
range of $m$ where the bound applies (with an appropriate larger constant $A_{m,n}$).
At each stage we construct a partitioning polynomial as before, but of a non-constant degree. We then use the bound
obtained in the previous stage to control the number of incidences inside of the cells of the polynomial partitioning. Finally, we use a separate argument (essentially the one given in the second part of the proof of Theorem \ref{th:arb}) to bound the number of
incidences with the points that lie on the zero set of the polynomial.
Each stage increases the constant of proportionality in the
bound by a constant factor, which is why the ``constant'' $A_{m,n}$ increases
as $m$ approaches $n^{3/2}$. For $j=1,2,\ldots,$
the $j$-th stage asserts the bound specified in the theorem when
$m\le n^{\alpha_j}$. The sequence of exponents $\{\alpha_j\}$
increases from stage to stage, and approaches $3/2$.
Each stage has its own constant of proportionality $A^{(j)}$.
The specific values of the exponents $\alpha_j$ (and the constants of proportionality) will be set later.
For the $0$-th, vacuous stage we use $\alpha_0=1/3$, and the bound
$O(n)$ that was noted above for $m\le n^{\alpha_0}$, with an implied initial constant of proportionality $A^{(0)}$.

In handling the $j$-th stage, we assume that
$n^{\alpha_{j-1}} < m \le n^{\alpha_j}$; if
$m\le n^{\alpha_{j-1}}$ there is nothing to do as we can use the (better) bound from the previous stage.
We construct an $r$-partitioning polynomial $f$, just as in the proof of Theorem \ref{th:arb}, except that its degree is not required to be a constant.
Put $\alpha=\alpha_{j-1}$.
To apply the bound from the previous stage uniformly within
each cell, we want to have a uniform bound on the number of circles
entering a cell. The average number circles entering a cell is proportional to $n/r^{2/3}$ (assuming
that the number of cells is $\Theta(r)$, an assumption made only for the sake of intuition).
A cell that intersects $tn/r^{2/3}$ circles, for $t>1$, induces $\lceil t\rceil$
subproblems, each involving all the points in the cell and up to
$n/r^{2/3}$ circles.
It is easily checked that the number of subproblems remains $O(r)$,
with a somewhat larger constant of proportionality, and that each subproblem
now involves at most $m/r$ points and at most $n/r^{2/3}$ circles.
Moreover, in cells that have strictly fewer than $n/r^{2/3}$ circles, we will assume that there are exactly $n/r^{2/3}$ circles, e.g., by adding dummy circles. This will not decrease the number of incidences.

We assume that the number of cells
is at most $br$, for some absolute constant $b$.
To apply the bound from the previous stage, we need to choose $r$ that will guarantee that
$$
\frac{m}{r} \le \left(\frac{n}{r^{2/3}}\right)^{\alpha} ,
\quad\quad\text{i.e.}\quad\quad
r^{1-2{\alpha}/3} \ge \frac{m}{n^{\alpha}} ,\quad\quad\text{i.e.}\quad\quad
r \ge \frac{m^{3/(3-2{\alpha})}}{n^{3{\alpha}/(3-2{\alpha})}} .
$$
We choose $r$ to be equal to the last expression. We note that
(i) $r\ge 1$, because $m$ is assumed to be greater than $n^\alpha$ and $\alpha < 3/2$, and
(ii) $r\le m$, because $m\le n^{3/2}$. Because of the somewhat weak bound
that we will derive below on the number of incidences with points that
lie on $Z(f)$ (the same bound as in the proof of Theorem~\ref{th:arb}),
this choice of $r$ will work only when $m$ is not
too large. The resulting constraint on $m$, of the form $m\le n^{\alpha_j}$,
will define the new range in which the bound derived in the present stage applies.

In more detail, the number of incidences within the partition cells is
\begin{align*}
I(\pts^\prime,\circles^\prime) &\le
A^{(j-1)} \sum_{i=1}^{br} \left( (m/r)^{3/7}(n/r^{2/3})^{6/7} +
  (m/r)^{2/3}(n/r^{2/3})^{1/2}q^{1/6} \right. \\
  & \qquad\qquad\qquad\qquad + \left. (m/r)^{6/11}(n/r^{2/3})^{15/22}q^{3/22}\log^{2/11}(m/r)
  + m/r + n/r^{2/3} \right) \\
& \le bA^{(j-1)} \left( m^{3/7}n^{6/7} +
  m^{2/3}n^{1/2}q^{1/6} + m^{6/11}n^{15/22}q^{3/22}\log^{2/11}m
  + m + nr^{1/3} \right) .
\end{align*}

We claim that our choice of $r$ ensures that
$nr^{1/3} \le m^{3/7}n^{6/7}$. That is,
$$
r^{1/3} =
\frac{m^{1/(3-2\alpha)}}{n^{\alpha/(3-2\alpha)}} \le
\frac{m^{3/7}}{n^{1/7}} .
$$
Indeed, this is easily seen to hold because $1/3 \le \alpha < 3/2$ and $m\le n^{3/2}$.
Recall that we also have $I(\pts_0,\circles^\prime) \le A'nr^{1/3}$ for some constant $A'$ (see (\ref{eq:c'p0})).
By choosing $A^{(0)}>A'$ (so that $A^{(j-1)}>A'$ for every $j$), we have
\[ I(\pts,\circles^\prime) = I(\pts_0,\circles^\prime) + I(\pts^\prime,\circles^\prime) \qquad \qquad \qquad \qquad \qquad \qquad \qquad \qquad \qquad \qquad \qquad \qquad \qquad\] \vspace{-13mm}

\begin{equation} \label{eq:stage}
\qquad \qquad \le 3bA^{(j-1)} \left( m^{3/7}n^{6/7} +
  m^{2/3}n^{1/2}q^{1/6} + m^{6/11}n^{15/22}q^{3/22}\log^{2/11}m
  + m \right).
\end{equation}
As proved in Theorem \ref{th:arb},
\begin{equation} \label{eq:stage2}
I(\pts_s,\circles_0) + I(\pts_p^{(2)},\circles_0) = O\left(mr^{2/3} + nr^{1/3}\right);
\end{equation}
this follows by substituting $D=O(r^{1/3})$ in the bounds in the proof of
Theorem \ref{th:arb}, which are $I(\pts_s,\circles_0)\le mD^{2}/2 + 2nD$ and
$I(\pts_p^{(2)},\circles_0) \le 44mD^2 + nD$.

It remains to bound $I(\pts_p^{(1)},\circles_0)$. For this, we again use an analysis similar to the one in Theorem \ref{th:arb}.
Let $f_1,f_2,\ldots,f_t$ be the factors of $f$ whose zero sets are
planes or spheres.
For $i=1,\ldots,t$, set $\pts_{p,i}^{(1)}=\pts_p \cap  Z(f_i)$  and
$m_{p,i}=|\pts_{p,i}^{(1)}|$. Let $n_{p,i}$ denote the number of circles of $\circles_0$ that are
fully contained in $Z(f_i)$ (ignoring, as before, circles that lie in
more than one of these surfaces).
Put $\pts_p^{(1)}=\bigcup_{i=1}^{t} \pts_{p,i}^{(1)}$.
Notice that (i) $t = O\left(r^{1/3}\right)$, (ii) $n_{p,i} \le q$
for every $i$, and (iii) $\sum_i n_{p,i} \le n$.
Applying \eqref{eq:oldBound}, we obtain
\begin{align}
I(\pts_p^{(1)},\circles_0) &=
\sum_{i=1}^t O\left( m_{p,i}^{2/3}n_{p,i}^{2/3} +
m_{p,i}^{6/11}n_{p,i}^{9/11}\log^{2/11}(m_{p,i}^3/n_{p,i})+m_{p,i}+n_{p,i}\right)  \nonumber \\
 &= \sum_{i=1}^t O\left(
 m_{p,i}^{2/3}n_{p,i}^{1/3}q^{1/3} +
 m_{p,i}^{6/11}n_{p,i}^{5/11}q^{4/11}\log^{2/11}(m_{p,i}^3)+m_{p,i}+n_{p,i}\right)  \nonumber \\
&= O\left( m^{2/3}n^{1/3}q^{1/3} +
m^{6/11}n^{5/11}q^{4/11}\log^{2/11}m+m+n \right), \label{eq:stage3}
\end{align}
where the last step uses H\"older's inequality.

We would like to combine \eqref{eq:stage}, \eqref{eq:stage2}, and \eqref{eq:stage3} to obtain the bound asserted in Theorem \ref{th:noEps}.
All the elements in these bounds add up to the latter bound, with an appropriate sufficiently large choice of $A^{(j)}$,
except for the term $O(mr^{2/3})$, which might exceed
the bound of the theorem if $m$ is too large. Thus, we restrict $m$ to satisfy
$$
mr^{2/3} \le m^{3/7}n^{6/7} , \quad\quad\text{i.e.}\quad\quad
r \le \frac{n^{9/7}}{m^{6/7}} .
$$
Substituting the chosen value of $r$, we thus require that
$$
\frac{m^{3/(3-2\alpha)}}{n^{3\alpha/(3-2\alpha)}}
\le \frac{n^{9/7}}{m^{6/7}} .
$$
That is, we require
$$
m \le n^{\frac{9+\alpha}{13-4\alpha}} .
$$
Recalling that we write the (upper bound) constraint on $m$
at the $j$-th stage as $m\le n^{\alpha_j}$, we have the recurrence
$$
\alpha_j = \frac{9+\alpha_{j-1}}{13-4\alpha_{j-1}} .
$$
To simplify this, we write $\alpha_j = \frac32 - \frac{1}{x_j}$, and
obtain the recurrence
$$
x_j = x_{j-1} + \frac47 ,
$$
with the initial value $x_0=\frac67$ (this gives the initial constraint
$m\le n^{1/3}$). In other words, we have $x_j = (4j+6)/7$, and
$$
\alpha_j = \frac32 - \frac{7}{4j+6} .
$$
The first few values are $\alpha_0 = 1/3$, $\alpha_1 = 4/5$,
$\alpha_2 = 1$, and $\alpha_3 = 10/9$.
Note that every $m<n^{3/2}$ is covered by the range of some stage.
Specifically, given such an $m$, it is covered by stage $j$, where
$j$ is the smallest integer satisfying
$$
m \le n^{\frac32 - \frac{7}{4j+6}} ,
$$
and straightforward calculations show that
$$
j = \left\lceil \frac32 \cdot \frac{\log(m/n^{1/3})}
{\log(n^{3/2}/m)} \right\rceil .
$$
Inspecting the preceding analysis, we see that the bound holds for the $j$-th
stage if we choose $A^{(j)}=A\cdot A^{(j-1)}$, where $A$ is a sufficiently
large absolute constant. Hence, for $m$ in the $j$-th range, the bound on
$I(\pts,\circles)$ has $A^{j+1}$ as the constant of proportionality.
This completes the description of the stage, and thus the proof of Theorem \ref{th:noEps}.
\end{theProof}

\noindent {\bf Remarks.} (1) The analysis holds for any $m<n^{3/2}$.
However, when $m$ is very close to $n^{3/2}$, say it is proportional to $n^{3/2}$,
then $j=\Theta(\log n)$, and the ``constant" $A^{(j)}$ is no longer a constant.
The requirement $m\le n^{3/2 -\delta}$ in the theorem is made to ensure that
$A_{m,n}$ does not exceed some constant threshold (which depends on $\delta$).
\medskip

\noindent (2) The case $m>n^{3/2}$ is not considered for this improvement, but we believe
that it too can be handled by similar techniques. Note that $m^{3/7}n^{6/7} = O(m)$
when $m>n^{3/2}$, so, ignoring the terms that depen on $q$, the overall bound is $O(m^{1+\eps})$,
for any $\eps>0$. One should be able to remove this dependence on $\eps$, as we did in the case
$m < n^{3/2}$.

\section{The number of popular points in an irreducible variety} \label{sec:reg}

It remains to prove Lemma \ref{le:reg}. To do so, we will use the three-dimensional
\emph{inversion transformation} $I:\RR^3 \to \RR^3$ about the origin (e.g.,
see \cite[Chapter 37]{Hart00}).  The transformation $I(\cdot)$ maps the point
$p=(x_1,x_2,x_3)\neq (0,0,0)$ to the point $\bar{p}=I(p)=(\bar{x}_1,\bar{x}_2,\bar{x}_3)$,
where
\[ \bar{x}_i = \frac{x_i}{x_1^2+x_2^2+x_3^2},\quad i=1,2,3. \]
A proof for the following lemma can be found in \cite[Chapter 37]{Hart00}.
\begin{lemma} \label{le:inv}
(a) Let $C$ be a circle incident to the origin. Then $I(C)$ is a line not passing through the origin. \\
(b) Let $C$ be a circle \emph{not} incident to the origin. Then $I(C)$ is a circle not passing through the origin. \\
(c) The converse statements of both (a) and (b) also hold. \hfill\(\Box\)
\end{lemma}

\begin{theProof}{of Lemma \ref{le:reg}}
Consider an irreducible surface $Z=Z(g)$ which is neither a plane nor a sphere,
and let $E=\operatorname{deg}(g)$.  Assume, for contradiction, that there exist
three popular points $z_1,z_2,z_3\in Z$. By translating the axes we may assume that
$z_1$ is the origin.  We apply the inversion transformation to $Z$. Since $I$ is
its own inverse, $I(Z)$ can be written as $Z(g\circ I)$. To turn $g \circ I$ into a
polynomial, we clear the denominators resulting from this transformation by multiplying
$g\circ I$ by a suitable (minimal) power of $x_1^2+x_2^2+x_3^2$. This does not
change the (real) zero set of $g\circ I$ (except for possibly adding the origin
$0$ to the set). We refer to the resulting polynomial as $\bar{g}$. Notice that the degree of
$\bar{g}$ is strictly smaller than $2E$, since to clear denominators we need
to multiply $g\circ I$ by at most $\left(x_1^2+x_2^2+x_3^2\right)^{E}$, and the
highest-degree terms will be the ones that were previously the linear terms (if they exist;
since $Z(g)$ contains the origin, $g$ has no constant term).
Since we have multiplied $g\circ I$ by the minimum power of $x_1^2+x_2^2+x_3^2$, we may
assume that $\bar g$ is not divisible by $x_1^2+x_2^2+x_3^2$.  If some other polynomial
divided $\bar g$, then after applying the inversion again and clearing denominators we
would obtain a non-trivial polynomial different from $g$ that divides $g$. Thus, since $g$ is irreducible, we
conclude that $\bar g$ is also irreducible.

By assumption, $Z(g)$ contains $44E^2$ circles incident to the origin
$z_1$. Lemma \ref{le:inv}(a) thus implies
that $Z\left(\bar g\right)$ contains at least $44E^2$ lines. We claim that
$Z(\bar g)$ is a ruled surface. If $\deg(\bar g)=2,$ then we can consider
all types of quadratic trivariate polynomials and observe that the ones
whose zero sets may contain more than $44E^2=176$ lines are all ruled
(namely, they are pairs of planes, cones, cylinders, 1-sheeted hyperboloids,
or hyperbolic paraboloids).
(If $\deg(\bar g)=1,$ $Z(\bar{g})$ is a plane.) If $\deg(\bar g)\ge3$, then
since $\bar g$ is irreducible of degree at most $2E$ and $Z(\bar g)$ contains
at least $44E^2 > 2E (11\cdot2E-24)$ lines, Corollary \ref{co:ruled} implies
that $Z(\bar g)$ is ruled. Since no regulus contains a point and more than
two lines through that point, $Z(\bar g)$ is not a regulus.
Moreover, $Z(\bar g)$ is not a plane since $Z(g)$ is neither a plane nor a sphere. That is, $Z(\bar g)$ is singly ruled.

Thus, $Z=Z(g)$ can be written as the union of a set of circles and a
(possibly empty) set of lines, all of which are incident to $z_1$;
these are the images under the inverse inversion of the lines spanning
$Z(\bar{g})$ (cf.~Lemma~\ref{le:inv}(c), and observe that
lines through $z_1$ are mapped to themselves by the inversion). By a
symmetric argument, this property also holds for $z_2$ and for $z_3$.
This implies that, for $i=1,2,3$, every point $u$ in $Z$ is incident
to a circle or a line that is also incident to $z_i$.
These three circles or lines are not necessarily distinct, but they can
all coincide only when $u$ lies on the unique circle or line $\gamma$ that
passes through $z_1,z_2,z_3$, and then all the above three circles or lines
coincide with $\gamma$.

The original surface $Z$
may or may not be ruled.  Recall that the only doubly ruled surfaces are
the hyperbolic paraboloid and the hyperboloid of one sheet.  Since both of
these surfaces do not contain a point that is incident to infinitely many
lines or circles contained in the surface, we conclude that $Z$ is not doubly
ruled. Since we have assumed that $Z$ is not a plane, it is not triply ruled
either.  Thus $Z$ is either not ruled or only singly ruled.

We define a point $u\in Z$ to be \emph{exceptional} if there are infinitely
many lines contained in $Z$ that pass through $u$ (think, e.g., of the case
where $Z(g)$ is a cone with this point as an apex). By Corollary 3.6 from
\cite{GK11}, if $Z$ is singly ruled, then $Z$ contains at most one
exceptional point. According to Corollary \ref{co:ruled}, if $Z$ is not
ruled, it contains only finitely many lines, and thus it cannot contain
any exceptional points. (Corollary~\ref{co:ruled}
does not apply to quadratic surfaces,
but it can be verified that the above property also holds in this case, by checking
all the possible types of quadratic surfaces.) Therefore, $Z$ contains
at most one exceptional point, and in particular we may assume that $z_2,z_3$
are not exceptional points. Since $z_2$ (resp., $z_3$) is popular but not
exceptional, there are infinitely many circles passing through $z_2$ (resp., $z_3$)
and contained in $Z$.  On the other hand, as already observed,
at most one circle can pass through the triplet
$z_1,z_2,z_3$.  Thus, after possibly interchanging the roles of $z_1, z_2$
and $z_3,$ we may assume that there exists an infinite collection of circles
contained in $Z$ that are incident to $z_2$ but not to $z_1$.

Consider the image $Z(\bar g)\subset\RR^3$ of $Z$ after applying
the inversion transform around the point $z_1$ (which we have translated
to become the origin) and let $\bar z_2 = I(z_2).$ According to
Lemma \ref{le:inv}(b), the infinite family of circles contained in $Z$
that are incident to $z_2$ but not to $z_1$ are transformed into an infinite
family of circles that are contained in $Z(\bar g)$ and incident to $\bar z_2$.
We denote the latter family as $\bar{\circles}$.

Consider a plane $\Pi$ and notice that $\Pi \cap Z(\bar{g})$ is an algebraic
curve of degree at most $2E$. This implies that $\Pi$ contains at most $E$
circles of $\bar{\circles}$.  Since this holds for any plane, there exists
an infinite subset $\bar{\circles}^\prime\subset \bar{\circles}$ such that
no two circles in $\bar{\circles}^\prime$ are coplanar.
($\Pi$ cannot be contained in $Z(\bar{g})$ since the latter surface is
irreducible, and if it is a plane then $Z$ is either a sphere or a plane.)

Since $Z(\bar g)$ is two-dimensional, by Theorem \ref{th:real} we have $(\bar g) = \mathbf{I}(Z(\bar g))$. Since $Z(\bar g)$ is ruled, Theorem \ref{th:flec2} implies $Z(\bar g) \subset Z(\fl{\bar g})\subset\RR^3$. Hence $\mathbf{I}(Z(\fl{\bar g})) \subset \mathbf{I}(Z(\bar g))$, and in particular $\fl{\bar g} \in \mathbf{I}(Z(\bar g)) = (\bar g)$. That is, $\bar g$ divides $\fl{\bar g}$ in $\RR[x_1,x_2,x_3]$, and thus also in $\CC[x_1,x_2,x_3]$. In particular, this means that $Z_{\CC}(\bar g)\subset Z_{\CC}(\fl{\bar g})$. By Theorem \ref{th:landsberg} this implies that $Z^*$ is also ruled.

In the remainder of the proof we will work mainly in complex projective
3-space $\CC\mathbf{P}^3$ instead of real affine space, which we have considered so far.

\paragraph{Projectivization.}
The \emph{projectivization} of a point $p=(p_1,p_2,p_3)\in\CC^3$ is
obtained by passing to \emph{homogeneous coordinates}, and by assigning
$p$ to $p^{\dagger}=(1,p_1,p_2,p_3)$ (where all nonzero scalar multiples
of $p^\dagger$ are identified with $p^\dagger$). To distinguish between
such homogeneous coordinates and coordinates in the affine spaces $\RR^3$
and $\CC^3$, we write them as $[x_0:x_1:x_2:x_3]$ (a rather standard notation,
see \cite{Mi95}), and will also use shortly a similar notation for the projective
Pl\"ucker 5-space of lines in 3-space.
The space of all points $[x_0:x_1:x_2:x_3]\neq 0$ is denoted as $\CC\mathbf{P}^3$.
As noted, two points $[x_0:x_1:x_2:x_3]$, $[x'_0:x'_1:x'_2:x'_3] \in \CC\mathbf{P}^3$
are considered to be equivalent if there exists a nonzero constant $\lambda\in \CC$
such that $x_0=\lambda x'_0$, $x_1=\lambda x'_1$, $x_2=\lambda x'_2$, and $x_3=\lambda x'_3$.
Given a point $[p_0:p_1:p_2:p_3]$ with $p_0\neq 0$, its \emph{dehomogenization}
with respect to $p_0$ is the affine point $(p_1/p_0,p_2/p_0,p_3/p_0)\in\CC^3$.
For more details, see \cite[Chapter 8]{CLO2}.

If $h\in\CC[x_1,x_2,x_3]$ is a polynomial of degree $E$, we can write
$h=\sum_{I} a_I x^I,$ where each index $I$ is of the form $(I_1,I_2,I_3)$
with $I_1+I_2+I_3\leq E$, and $x^I=x_1^{I_1}x_2^{I_2}x_3^{I_3}.$ Define
\begin{equation*}
h^{\dagger} = \sum_I a_I x_0^{E-I_1-I_2-I_3}x_1^{I_1}x_2^{I_2}x_3^{I_3}.
\end{equation*}
Then $h^{\dagger}$ is a homogeneous polynomial of degree $E$, referred to
as the \emph{homogenization} of $h$. We define the \emph{projectivization}
of the complex surface $Z(h)$ to be the zero set of $h^\dagger$ in the
three-dimensional complex projective space $\CC\mathbf{P}^3$.
We define the \emph{complex projectivization} of a real surface $S=Z(h)$
to be the projectivization of the complexification $S^*$ of $S$.

Let $\hat{Z}\subset \CC\mathbf{P}^3$ be the complex projectivization
of the surface $Z(\bar g)$.

For the next steps of the analysis, we introduce the so-called
\emph{absolute conic} in $\CC\mathbf{P}^3$ (e.g., see \cite{NS12})
\[
\Gamma = \bigg\{[x_0:x_1:x_2:x_3]\ \Big\vert \ x_0=0,\ x_1^2+x_2^2+x_3^2=0 \bigg\} .
\]
Notice that $\Gamma$ is contained in the plane
at infinity $x_0=0$, and does not contain any real point (a point all of
whose coordinates are real).

We will need the following simple lemma.
\begin{lemma} \label{le:derivatives}
Let $f\in \CC[x_0,x_1,x_2,x_3]$ be a homogeneous polynomial of degree $D$,
and let $S = Z(f) \subset \CC\mathbf{P}^3$. Let $p$ be a point in $S$, let
$v$ be a direction in $\CC\mathbf{P}^3$, and let $\ell$ be the line incident
to $p$ with direction $v$. If all partial derivatives of $f$ of order at most
$D$ vanish at $p$ in the direction $v$, then $\ell \subset S$.
\end{lemma}
\begin{theProof}{\!\!}
Let $f_0$ be the restriction of $f$ to the line $\ell$, and notice that
$f_0$ is a univariate polynomial of degree at most $D$.
Since all derivatives of $f_0$ of degree
at most $D$ vanish at $p$, $f_0$ must be identically 0.
\end{theProof}

\begin{lemma} \label{le:irredRuled}
$\hat{Z}\subset \CC\mathbf{P}^3$ is irreducible and singly ruled.
\end{lemma}
\begin{theProof}{\!\!}
According to the above discussion, the complexification $Z^*$ of $Z(\bar g)$
is irreducible and ruled. It is singly ruled since otherwise $Z^*$ would
be a complex plane or regulus,
and it cannot be that the real part of a complex plane or regulus is a singly
ruled surface.
It remains to consider the projectivization $Z^*$ of the complexification.

If the homogenization $f^\dagger$ of a non-constant polynomial $f$ is
divisible by a polynomial $f_0$, then $f$ is divisible by the polynomial
obtained by substituting $x_0=1$ in $f_0$. Moreover, $f^\dagger$ is not
divisible by any polynomial of the form $x_0^a$, for an integer $a>0$
(since $f$ is not a constant).  Thus, if $f$ is irreducible, then so is
$f^\dagger$. This in turn implies that $\hat{Z}$ is irreducible.

By definition, for every point $p$ in the complexification $Z^*$ of $Z(\bar g)$
there exists a line $\ell$ that is incident to $p$ and fully contained
in $Z^*$. Notice that the projective point $p^\dagger\in\hat{Z}$ is incident
to the projective line $\ell^\dagger\subset \hat{Z}$ (i.e., the locus of the
projectivizations of the points of $\ell$). Thus, for any point
$p=[p_0:p_1:p_2:p_3] \in\hat{Z}$ for which $p_0\neq 0$, there is a line
incident to $p$ and fully contained in $\hat{Z}$.  It remains to consider
points $p=[p_0:p_1:p_2:p_3] \in\hat{Z}$ for which $p_0 = 0$.
Since $\hat Z$ is irreducible and cannot be the plane $Z(x_0)$, then $\hat Z \cap Z(x_0)$ is a 1--dimensional curve. Since $p$ must have at least one non-zero coordinate, we may assume without loss of generality that $p_3\neq 0$. Let $Z_0\subset\CC^3$ be the dehomogenization $\hat Z\backslash \{x_3=0\}$ with respect to $x_3$, and note that the image of p is contained in $Z_0$. Notice that every projective line that is fully contained in $Z$ but not contained in $Z(x_3)$ corresponds to a line fully contained in $Z_0$. Since $\hat Z$ contains a line through every point with a nonzero $x_0$-coordinate (and the set of points with zero $x_0$-coordinate form a proper sub-variety of $\hat Z$), $Z_0$ contains infinitely many lines. Applying Theorem 2.12 implies that $Z_0$ is ruled. Thus, $Z_0$ contains a line $\ell_0$ incident to $p$. The projectivization of $\ell_0$ is fully contained in $\hat Z$ and incident to $p$. In conclusion, $\hat Z$ is ruled; it remains to show that it is singly ruled.

Notice that there is a bijection between the lines in affine space
and the lines in projective space that are not fully contained in the
plane $Z(x_0)$.
Consider a projective
point $p^\dagger \in Z(\bar g^\dagger) \setminus Z(x_0)$ and its corresponding
dehomogenized affine point $p$. If $p^\dagger$ is incident to two lines that
are fully contained in $Z(\bar g^\dagger)$, then the two corresponding affine lines are incident
to $p$ and are fully contained in $Z(\bar g)$.
Thus, $Z(\bar g^\dagger)$ is singly ruled, since if $Z(\bar g^\dagger)$
were doubly or triply ruled then $Z(\bar g)$ would also have to be doubly
or triply ruled, which is not the case.
\end{theProof}

\paragraph{Using the Pl\"ucker representation of lines.}
With all these preparations, we reach the following scenario.
We have an irreducible singly-ruled surface $\hat Z$ in $\CC\mathbf{P}^3$
(the complex projectivization of $Z(\bar g)$), which contains an infinite
family $\hat{C}$ of circles (the complex projectivization of the circles
of $\bar{C}$), no pair of which are coplanar.
The following arguments are based on the recent
work of Nilov and Skopenkov \cite{NS12} concerning surfaces that
are ``ruled" by lines and circles.  Before taking the circles into
account, we first probe deeper into the structure of our ruled surface
$\hat Z$ by considering the Pl\"ucker representation of lines in 3-space.
Specifically, let
\begin{equation}\label{representationOfPluckerQuadric}
\Lambda=\biggl\{[x_0:\ldots:x_5] \mid x_0x_5+x_1x_4+x_2x_3=0\biggr\}\subset\CC\mathbf{P}^5
\end{equation}
be the \emph{Pl\"ucker quadric}.
Given a point $p=[x_0:\ldots:x_5]\in\Lambda$, at least two of the four
``canonical" points $[0:x_0:x_1:x_2]$, $[x_0:0:-x_3:-x_4]$,
$[x_1:x_3:0:-x_5]$, and $[x_2:x_4:x_5:0]$ cannot be the undefined point
$[0:0:0:0]$ because each of the six coordinates $x_0,\ldots,x_5$, not all zero,
appears as a coordinate of two of these points.
Then there exists a unique line $\ell_p$ in
$\CC\mathbf{P}^3$ that passes through all nonzero canonical points of $p$.
We refer to the map $p\to \ell_p$ as the \emph{Pl\"ucker map}, and observe
that it is a bijection between the points $p\in \Lambda$ and the lines
$\ell_p \subset \CC\mathbf{P}^3$.  Further details about the Pl\"ucker map
and the Pl\"ucker quadric can be found, e.g., in \cite[Section 8.6]{CLO2}.

Let $\Lambda_{\hat Z}=\{p\in \Lambda \mid \ell_p\subset \hat Z\}$;
that is, $\Lambda_{\hat Z}$ is the set of all points in $\Lambda$ that
correspond to lines that are fully contained in $\hat Z$. We claim that
$\Lambda_{\hat Z}$ is an algebraic variety in $\CC\mathbf{P}^5$ that is
composed of a single one-dimensional irreducible component, possibly
together with an additional finite set of points. Example 6.19 of \cite{harris} establishes that $\Lambda_{\hat Z}$ is a projective variety (indeed, it is an example of a Fano variety). We must now show that an irreducible two-dimensional surface in $\CC\mathbf{P}^3$
that does not contain any planes, cannot contain a two-dimensional family of lines.
This implies that $\Lambda_{\hat Z}$ is a one-dimensional set.

\begin{lemma} \label{le:pluckerDimension}
Let $S \subset \CC\mathbf{P}^3$ be a ruled surface that does not contain
any planes and let $\gamma \subset \CC\mathbf{P}^5$ be the set of points
on the Pl\"ucker quadric that correspond to lines contained in $S$.
Then $\gamma$ is one-dimensional.
\end{lemma}
\begin{theProof}{\!\!}
Assume, for contradiction, that $\gamma$ contains a two-dimensional
irreducible component $\gamma_2$, and let $\Pi \subset \CC\mathbf{P}^3$
be a generic plane. Since we are in projective space, $\Pi$ intersects
every line of $\gamma_2$.  Each line contained in $S$ whose corresponding
point is in $\gamma_2$ intersects the curve $\sigma = S \cap \Pi$.
Since the singular points of $S$ are contained in a one-dimensional
variety (this follows, e.g., from Sard's lemma \cite{Sard42}), and
since $\Pi$ is a generic plane, we may assume that $\Pi$ contains only
finitely many singular points of $S$.  Let $F: \gamma_2 \to \sigma$ be
a mapping that sends each point $p\in \gamma_2$ to the intersection point of
$\sigma$ with the line corresponding to $p$.

We claim that $\sigma$ contains
a non-singular point $q$ of $S$ such that $F^{-1}(q)$ is infinite (i.e.,
infinitely many lines that correspond to points of $\gamma_2$ are incident
to $q$). Indeed, first notice that it is impossible for $F^{-1}(q)$ to be two-dimensional, for any $q\in \sigma$.
Indeed, if $F^{-1}(q)$ were two-dimensional, then the intersection of the corresponding lines with any sphere around
$q$ would also be two-dimensional, contradicting the fact that $S$ is two-dimensional.
Therefore, for every point $q\in \sigma$, $F^{-1}(q)$ is either one-dimensional or zero-dimensional.
Moreover, it cannot be that all of these preimages are zero-dimensional, since a one-dimensional set of
zero-dimensional varieties cannot cover the entire two-dimensional family of lines;
see, e.g., \cite[Exercise 3.22(b)]{Hart83}.
Thus, there are infinitely many points $q\in \sigma$, such that $F^{-1}(q)$ is infinite.
Since there are finitely many singular points of $S$ in $\sigma$,
there are non-singular points with such an infinite preimage.

We can now complete the proof of Lemma \ref{le:pluckerDimension}. We have a smooth point $q$ of $S$ such that there are infinitely many lines that pass through $q$ and are contained in $S$. These lines must also be contained in the tangent plane $T_qS$. Thus by Lemma \ref{le:commonLines}, $S$ must contain the tangent plane $T_qS$. This contradicts
the fact that $S$ does not contain any planes.

The above proof still holds when assuming that the dimension of $\gamma$
is larger than 2.  Notice that $\gamma$ cannot be zero-dimensional, since then
the higher-dimensional extension of B\'ezout's theorem (Theorem \ref{th:bezout3})
would imply that it is finite, and the union of the corresponding lines will not be two-dimensional
\end{theProof}

Lemma \ref{le:pluckerDimension} implies that $\Lambda_{\hat Z}$ is one-dimensional. To prove that it consists of a
single one-dimensional irreducible component (possibly with additional zero-dimensional components), we shall first require the following to Lemmas.

\begin{lemma} \label{le:linesUnion}
Let $\gamma \subset \CC\mathbf{P}^5$ be a projective subvariety of the Pl\"ucker quadric $\Lambda$.
Then $\bigcup_{p\in \gamma}\ell_p \subset \CC\mathbf{P}^3$ is also a projective variety.
\end{lemma}
Lemma \ref{le:linesUnion} is a special case of Proposition 6.13 from \cite{harris}

\begin{lemma} \label{le:crossEll}
Let $\ell \subset\CC\mathbf{P}^3$ be a projective line, and let $V\subset \CC\mathbf{P}^5$ be the set of points on
the Pl\"ucker quadric  $\Lambda$ that correspond to lines that intersect $\ell$. Then $V$ is a projective variety.
\end{lemma}
Lemma \ref{le:crossEll} is a special case of Example 6.14 from \cite{harris}.

Let $\gamma$ be an irreducible one-dimensional component of $\Lambda_{\hat Z}$.
According to Lemma \ref{le:linesUnion}, $\bigcup_{p\in \gamma}\ell_p$ is a two-dimensional
algebraic variety that is fully contained in $\hat Z$.  Since $\hat Z$ is
irreducible, $\bigcup_{p\in \gamma}\ell_p=\hat Z$ and $\gamma$ corresponds
to a generating family of $\hat Z$.  Let $q\in \CC\mathbf{P}^{5}$ be a point
of $\Lambda_{\hat Z} \setminus \gamma$, and let $\ell_q \subset \CC\mathbf{P}^{3}$
be the line that corresponds to $q$.  According to Lemma \ref{le:crossEll}, the
set $V_q \subset \CC\mathbf{P}^{5}$, of points that correspond to lines that
intersect $\ell_q$, is a variety.  Note that $\gamma \cap V_q$ is infinite,
because every point on $\ell_q$ lies on some generator line $\ell_p$ for
$p\in \gamma$. Since $\gamma$ is irreducible, then $\gamma \subset V_q$.
That is, every line in the generating family $ \{\ell_p\}_{p\in\gamma}$ of
${\hat Z}$ intersects $\ell_q$. If there are at least three points in
$\Lambda_{\hat Z} \setminus \gamma$, then each line in the generating
family of ${\hat Z}$ intersects three given lines, which implies that
${\hat Z}$ is either a regulus or a plane.\footnote{A nice proof for this claim, which holds in $\RR^3, \CC^3$, and $\CC\mathbf{P}^{3}$, can be found in \url{http://math.mit.edu/~lguth/PolyMethod/lect10.pdf} (version of June 2013).}
Since reguli and planes are not singly ruled, it follows that
$\Lambda_{\hat Z}$ is composed of an irreducible one-dimensional curve,
and at most two other points (the additional points correspond to
non-generating lines that are fully contained in $\hat Z$).\footnote{%
This also implies that the reguli are the only doubly ruled surfaces in $\CC\mathbf{P}^{3}$.}

\paragraph{Adding $\Gamma$ to the analysis.}
Consider a line $\ell$ whose pre-image under the Pl\"ucker map is the
point $[p_0:\ldots:p_5]\in \Lambda$, such that $\ell$ does not lie in the plane
at infinity. Then $\ell$ intersects $\Gamma$ if and only if
$p_0^2+p_1^2+p_2^2=0$. Indeed, recall that $\ell$ contains the point
$[0:p_0:p_1:p_2]$, and this is the only point of $\ell$
on the plane at infinity $Z(x_0)$, for otherwise $\ell$ would be fully
contained in that plane. (It cannot be that $p_0=p_1=p_2=0$, since then all
four points $[0:p_0:p_1:p_2]$, $[p_0:0:-p_3:-p_4]$, $[p_1:p_3:0:-p_5]$, and
$[p_2:p_4:p_5:0]$ would have a zero $x_0$-coordinate, implying that $\ell$ is
contained in $Z(x_0)$.) Thus, the set
$\Gamma_{\Lambda}=\{p\in\Lambda \mid \ell_p\cap\Gamma\neq\emptyset\}$ is
an algebraic variety of codimension $1$ in $\Lambda$.
Since the irreducible one-dimensional component of $\Lambda_{\hat Z}$ is also
a variety, either it is fully contained in $\Gamma_{\Lambda}$, or the intersection
$\Lambda_{\hat Z}\cap \Gamma_{\Lambda}$ is a zero-dimensional variety, and
therefore finite according to the higher-dimensional extension of B\'ezout's theorem (Theorem \ref{th:bezout3}).
If the former case occurs, then at most two lines in $\hat Z$
do not intersect $\Gamma$. However, since $\hat Z$ is the complex
projectivization of a real ruled surface, $\hat Z$ contains
infinitely many real lines (lines whose defining equations involve
only real coefficients) that are not contained in the plane $Z(x_0)$, and
if $\ell$ is such a line then
$\ell\cap Z(x_0)$ consists of real points. This is a contradiction
since the curve $\Gamma$ contains no real points. Therefore,
the intersection $\Lambda_{\hat Z}\cap \Gamma_{\Lambda}$ is finite.

\begin{lemma} \label{le:conic}
Every line intersects $\Gamma$ in at most two points.
\end{lemma}
\begin{theProof}{\!\!}
By B\'ezout's theorem (Theorem \ref{th:bezout}), applied in the plane at
infinity $h=Z(x_0)$, $\Gamma$ has at most two intersection points with any
line that is not fully contained in $\Gamma$ (clearly, lines not contained in $h$
can meet $\Gamma$ at most once). Thus, it suffices to prove that no line is
fully contained in $\Gamma$.

Regard $h$ as the standard complex projective plane $\CC\mathbf{P}^2$, with
homogeneous coordinates $[x_1:x_2:x_3]$. A line $\ell\in h$ has an equation of
the form $a_1x_1+a_2x_2+a_3x_3 = 0$, and there exists at least one coordinate,
say $x_3$, with $a_3\neq 0$. This allows us to write the equation of $\ell$ as
$x_3 = \alpha x_1 + \beta x_2$, so its intersection with $\Gamma$ satisfies the equation
\[ x_1^2+x_2^2+(\alpha x_1 + \beta x_2)^2= 0, \quad \text{or} \quad
(1+\alpha^2)x_1^2 + 2\alpha\beta x_1 x_2 + (1+\beta^2)x_2^2 = 0. \]
This is a quadratic equation, whose coefficients cannot all vanish, as is
easily checked. Hence it has at most two solutions, which is what the lemma asserts.
\end{theProof}

Lemma \ref{le:conic} implies that $\Gamma\cap\hat Z$ is a finite set.
Indeed, if this
were not the case, then there would exist infinitely many points of $\Gamma$
that lie in $\hat{Z}$ and each of them is therefore incident to a line contained
in $\hat Z$. Since every line meets $\Gamma$ in at most two points, $\Gamma$
would have intersected infinitely many lines contained in $\hat{Z}$. This is a
contradiction since, as argued above, $\Lambda_{\hat Z}\cap \Gamma_{\Lambda}$
is a finite intersection.

\paragraph{Adding the circles to the analysis.}
Let $\bar{\mathcal C}^\prime$ be the collection of circles described earlier;
that is, an infinite set of pairwise non-coplanar circles that are fully
contained in $Z(\bar g)$ and incident to $\bar{z}_2$.
Let $\hat{\mathcal C}^\prime$ be the corresponding collection of
the complex projectivizations of these circles. As just argued,
all of the intersection points between the circles of
$\hat{\mathcal C}^\prime$ and $\Gamma$ must lie in the finite
intersection $\Gamma\cap\hat Z$.
\begin{lemma} \label{le:circlesGamma}
Each circle $\hat{C}$ in
$\hat{\mathcal C}^\prime$ intersects $\Gamma$ in
precisely two points.
\end{lemma}
\begin{theProof}{\!\!}
Each circle $\hat{C}$ in
$\hat{\mathcal C}^\prime$ is the complex projectivization of a real
circle $C$.
Consider $C$ as the intersection of its supporting plane, whose equation
is $ax_0+bx_1+cx_2+dx_3=0$, for appropriate parameters $a,b,c,d$, with
some suitable sphere whose equation is
$(x_1-a'x_0)^2 + (x_2-b'x_0)^2 + (x_3-c'x_0)^2=d'x_0^2$, for appropriate
parameters $a',b',c',d'$.
Since $C$ is real, all coefficients of these equations can be assumed
to be real, and $d'>0$.

By combining these equations of $C$ with the equations $x_0=0$ and
$x_1^2+x_2^2+x_3^2=0$ of the absolute conic $\Gamma$, we obtain the system
$bx_1+cx_2+dx_3=0$, $x_1^2+x_2^2+x_3^2=0$ (where the second equation arises
twice), which always has two distinct solutions when $b,c,d$ are real.
Indeed, using the notation in the proof of Lemma \ref{le:conic},
the intersection points satisfy a quadratic equation of the form
$(1+\alpha^2)x_1^2 + 2\alpha\beta x_1 x_2 + (1+\beta^2)x_2^2 = 0$
(or a similar, symmetrically defined equation in another pair of variables)
where $\alpha$ and $\beta$ are real. The discriminant of this equation is
$4\alpha^2\beta^2-4(1+\alpha^2)(1+\beta^2) = -4(\alpha^2+\beta^2+1)$, which
is always nonzero (and strictly negative) when $\alpha$ and $\beta$ are real.
Also, the coefficients of $x_1^2,x_2^2$ are both nonzero, and thus the
equation has exactly two (complex conjugate) solutions.
\end{theProof}
\paragraph{The final stretch.}
Since $\hat C^\prime$ contains infinitely many circles and  $\Gamma
\cap \hat Z$ is finite, by the pigeonhole principle there must exist two circles
$C_1,C_2$ in $\hat{\mathcal C}^\prime$ such that the sets $C_1\cap \Gamma$
and $C_2\cap\Gamma$ are identical (each being a set of two points).
By construction, $C_1$ and $C_2$ are contained in two
distinct planes $\Pi_1$ and $\Pi_2$. Consider the line
$\ell= \Pi_1 \cap \Pi_2$ and notice that it contains $C_1 \cap C_2$.
Thus, $\ell$ contains the two intersection points of $C_1,C_2$ with
$\Gamma$. Since these two points are contained in the plane $\{x_0=0\}$,
$\ell$ is also contained in this plane. However, this is impossible,
since $\ell$ also contains $\bar{z}_2$ (common to all circles of
$\hat{\mathcal C}^\prime$), which is not in the plane $\{x_0=0\}$.
This contradiction completes, at long last, the proof of Lemma \ref{le:reg}.
\end{theProof}

\section{Unit circles} \label{sec:uni}

In this section we consider the special case where all circles in
$\circles$ have the same radius, say $1$. The analysis is very similar
to the general case, except for two key issues: (a) In the general case
we have used the fact that the incidence graph in $\pts\times\circles$
does not contain $K_{3,2}$ as a subgraph, to derive the weaker
``bootstrapping'' bound $I(\pts,\circles) = O(n^{2/3}m+n)$.
Here, in Lemma~\ref{le:53} below, we replace this estimate by an
improved one, exploiting the fact that all circles are congruent.\footnote{%
  If all our circles were coplanar or cospherical, life would have been
  simpler, since then the incidence graph does not contain $K_{2,3}$
  as a subgraph, which is the basis for deriving the improved planar
  bound $I(\pts,\circles) = O(m^{2/3}n^{2/3}+m+n)$. In three dimensions
  the incidence graph can contain $K_{2,q}$ for any value of $q$,
  making the analysis more involved and subtler.}
(b) When considering the case of circles that lie in a common
plane or sphere, we use the improved planar bound for unit circles
$I(\pts,\circles) = O(|\pts|^{2/3}|\circles|^{2/3}+|\pts|+|\circles|)$
(e.g., see~\cite{Sz}). These two improvements result in the
sharper bound of Theorem~\ref{th:uni}, which we restate here for
the convenience of the reader.

\paragraph{Theorem \ref{th:uni}}
\emph{Let $\pts$ be a set of $m$ points and let $\circles$ be a set of $n$
unit circles in $\RR^3$, let $\eps$ be an arbitrarily small positive constant,
and let $q\leq n$ be an integer.  If no plane or sphere contains more than $q$
circles of $\circles$, then
\[
I(\pts,\circles) = O\left(m^{5/11+\eps}n^{9/11} +
m^{2/3+\eps}n^{1/2}q^{1/6} + m + n\right),
\]
where the constant of proportionality depends on $\eps$.
}\vspace{3mm}

\hspace{-10.5mm} \begin{theProof}{\!\!}
We first establish the following lemma, which improves
the weaker bound on $I(\pts,\circles)$, as discussed in (a) above.
\begin{lemma} \label{le:53}
Let $\pts$ be a set of $m$ points in $\RR^3$ and let
$\circles$ be a collection of unit circles in $\RR^3$, so that
each circle of $\circles$ is incident
to at least three points of $\pts$. Then $|\circles| = O(m^{5/2})$.
A stronger statement is that the number of circles of $\circles$ that pass
through any fixed point $o\in \pts$ and through at least two other
points is $O(m^{3/2})$.
\end{lemma}

\begin{theProof}{\!\!}
It clearly suffices to establish only the second claim of the lemma.
Fix one point $o$ of $\pts$ and let $\pts^\prime = \pts\backslash\{o\}$. For each $a\in\pts^\prime,$ let $\sigma_a$ be the locus of all points $w\in\RR^3$ such that $o$, $a$, and $w$ lie on a common unit circle. The set $\sigma_a$ is an
algebraic surface of revolution, obtained by taking any unit circle
passing through $o$ and $a$ and by rotating it around the line $oa$.
If $o$ and $a$ are diametral, that is, if $|oa|=2$, then $\sigma_a$ is a
sphere. If $|oa|>2$ then $\sigma_a$ is empty. Otherwise, $\sigma_a$
is easily seen to be an irreducible surface of degree $4$; the ``outside''
portion of $\sigma_a$ resembles a sphere pinched at $o$ and $a$, which
are the only singular points of $\sigma_a$;
the ``inner" portion resembles a pointy (American) football.

Let $\mathcal S = \{\sigma_a\colon a\in\pts^\prime\}$. In order to prove the second claim of the Lemma, it suffices to show that $I(\pts^\prime,\mathcal{S})=O(m^{3/2})$. We require with the following lemma.
\begin{lemma}
There exists an absolute constant $s$ such that for all triples $a,b,c\in\pts^\prime$, we have $|\sigma_a\cap\sigma_b\cap\sigma_c|\leq s$
\end{lemma}
\begin{theProof}{\!\!}
If $|\sigma_a\cap\sigma_b\cap\sigma_c|$ is finite, then  Milnor's
theorem (Theorem \ref{th:Milnor}) implies that this number is at most
some constant $E$. By setting $s$ to be, say, $E+1$, we ensure that the
intersection must be infinite.

\begin{figure}[h]
\centerline{\placefig{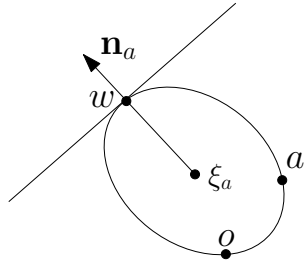}{0.26\textwidth}}
\vspace{-1mm}

\caption{\small \sf The normal $\nn_a$ to $\sigma_a$ lies on the ray
$\vec{\xi_aw}$, where $\xi_a$ is the center of the unit circle passing
through $o$, $a$, and $w$.}
\label{fi:lim}
\vspace{-2mm}
\end{figure}

Consider then the case where $\sigma_a \cap \sigma_b \cap \sigma_c$ is
a one-dimensional curve $\gamma$ (it cannot be two-dimensional because
$\sigma_a$, $\sigma_b$, $\sigma_c$ are distinct irreducible varieties,
no pair of which can overlap in a two-dimensional subset),
and let $w$ be a smooth point on $\gamma$.
Let $\tau$ be the tangent to $\gamma$ at $w$.
Then $\tau$ is orthogonal to the three respective normals $\nn_a$,
$\nn_b$, $\nn_c$ to $\sigma_a$, $\sigma_b$, $\sigma_c$ at $w$. In
other words, these normals must be coplanar. Now, because $\sigma_a$
is the surface of revolution of a circle, $\nn_a$ lies on the ray
$\vec{\xi_aw}$, where $\xi_a$ is the center of the unit circle passing
through $o$, $a$, and $w$; an illustration is provided in Figure \ref{fi:lim}.
Symmetric properties hold for $\sigma_b$
and $\sigma_c$, with respective centers $\xi_b$, $\xi_c$.

In other words, the argument implies that $w$, $\xi_a$, $\xi_b$,
and $\xi_c$ all lie in a common plane $h$. However, all three
centers $\xi_a$, $\xi_b$, and $\xi_c$ must lie on the perpendicular
bisector plane $\pi$ of $ow$, which does not contain $w$, so
$\pi\ne h$, and these centers then
have to lie on the intersection line $\ell=h\cap \pi$.
This is impossible if $\xi_a$, $\xi_b$, and $\xi_c$ are distinct,
because it is impossible for three distinct collinear points to
be at the same distance (namely, $1$) from $o$. Assume then that
$\xi_a=\xi_b$, say. That is, we have two distinct unit circles passing
through $o$ and $w$ with a common center $\xi=\xi_a=\xi_b$, which is
possible only when $|ow|=2$ (that is, $ow$ is a diameter of both circles).
Moreover, $\xi$ lies at distance $1$ from $o$, $a$, and $b$, so it is
the center of a unit ball that passes through these points. There
can be at most two such balls, so there are only two possible locations
for $\xi$. Since $\xi$ is the midpoint of $ow$ (recall that $ow$ is a
diameter of the sphere $\sigma_a$), it follows that
there are only two possible locations for $w$.
That is, $\gamma$ has at most two smooth points,
which is impossible, as follows, say, from Sard's theorem (e.g., see~\cite{Sard42}).
\end{theProof}

We can now apply Theorem 2 from \cite{Zahl11} to conclude that $I(\pts^\prime,\mathcal{S})\leq |\pts^\prime|^{3/4}|\mathcal{S}|^{3/4}+|\pts^\prime|+|\mathcal{S}|$. Since $|\pts^\prime|=O(m)$ and $|\mathcal{S}|=O(m)$, Lemma \ref{le:53} follows.
\end{theProof}

Consider pairs of the form $(p,c)$ where $p \in \pts$ and the circle $c\in \circles$ is incident to $p$ and to at least two other points of $\pts$.
By Lemma~\ref{le:53}, every point of $\pts$ can participate in at most $O(m^{3/2})$ such pairs, and thus the number of pairs is $O(m^{5/2})$.
This implies that $I(\pts,\circles) = O(m^{5/2}+n)$, so it is $O(n)$ for $m=O(n^{2/5})$
(recall that in the general case this could be claimed only for $m=O(n^{1/3})$).

The proof of Theorem~\ref{th:uni} now proceeds in complete analogy
with the proof of Theorem~\ref{th:arb}, except for the modifications
mentioned in (a) and (b) above. Specifically, we construct an
$r$-partitioning polynomial, of degree $O(r^{1/3})$, for a
sufficiently large constant parameter $r$, and consider
separately points of $\pts$ in the cells of the partition,
and points on $Z(f)$. The bound for the former kind of points
is handled via induction, in much the same way as before, except
that we replace the term $O(n)$, towards the derivation of
(a bound analogous to the one in) (\ref{eq:c'p'2}),
by $O(m^{5/11}n^{9/11})$, which holds for $n=O(m^{5/2})$.
We also remove the terms of the form $O(m^{6/11+\eps}n^{9/11})$
(see below for a justification). This results in the modified bound
\begin{equation} \label{eq:c'p'2-unit}
I(\pts_0 \cup \pts^\prime,\circles^\prime) \le
\frac{\alpha_1}{3}\left(m^{5/11+\eps}n^{9/11} +
m^{2/3+\eps}n^{1/2}q^{1/6}\right) +\alpha_2 |\pts'|.
\end{equation}
This ``explains'' why we can use here the improved exponents
$5/11$ and $9/11$ instead of the weaker respective ones $3/7$ and $6/7$.

The second modification is in handling incidences involving private
points on $Z(f)$ that lie in planes or spheres that are zero sets of
respective irreducible factors of $f$. Here, in the derivation of
a bound analogous to the one in (\ref{eq:reg1}), we use the sharper planar bound
$I(\pts,\circles) = O(|\pts|^{2/3}|\circles|^{2/3}+|\pts|+|\circles|)$,
which also holds when the points and circles are all cospherical.
This replaces (\ref{eq:reg1}) by the sharper bound
\begin{equation}\label{eq:reg1-unit}
I(\pts_p^{(1)},\circles_0) \le
\frac{\alpha_1}{3} m^{2/3+\eps}n^{1/2}q^{1/6} +
\alpha_2(m_p^{(1)}+n/3).
\end{equation}
The rest of the analysis remains unchanged, and leads to the bound
asserted in the theorem.
\end{theProof}

By applying a the techniques presented in Section \ref{ssec:eps} we obtain the following theorem.

\begin{theorem} \label{th:UnitNoEps}
Let $\pts$ be a set of $m$ points and let $\circles$ be a set of $n$
unit circles in $\RR^3$, let $q\leq n$ be an integer, and let $m=O(n^{3/2-\delta})$, for
some fixed arbitrarily small constant $\delta>0$.
If no sphere or plane contains more than $q$ circles of $\circles$, then
\begin{equation*}
I(\pts,\circles) \le A_{m,n}\left(m^{5/11}n^{9/11} +
m^{2/3}n^{1/2}q^{1/6} + m + n\right),
\end{equation*}
where $A_{m,n} = A^{\left\lceil \frac{\log{(m^5/n^2)} }{3\log{(n^{3/2}/m)} } \right\rceil+1}$,
for some absolute constant $A>1$. \hfill\(\Box\)
\end{theorem}

Since the proof of Theorem \ref{th:UnitNoEps} is almost identical to the proof of Theorem \ref{th:noEps}, we omit it.

\section{Applications}\label{furtherApplicationsSection}

\paragraph{High-multiplicity points.}
The following is an easy but interesting consequence of
Theorems \ref{th:arb} and \ref{th:uni}.
\begin{corollary}
(a) Let $\circles$ be a set of $n$ circles in $\RR^3$, and let $q\leq n$ be an
integer so that no sphere or plane contains more than $q$ circles of
$\circles.$ Then there exists a constant $k_0$
(independent of $\circles$) such that for any $k\ge k_0$,
the number of points incident to at least $k$ circles of $\circles$ is
\begin{equation} \label{eqnew}
\tilde O\left(\frac{n^{3/2}}{k^{7/4}} + \frac{n^{3/2}q^{1/2}}{k^3} +
\frac{n^{3/2}q^{3/10}}{k^{11/5}} + \frac{n}{k}   \right).
\end{equation}
In particular, if $q=O(1)$, the number of such points is
\[ \tilde O\left(\frac{n^{3/2}}{k^{7/4}} + \frac{n}{k}   \right). \]
(b) If the circles of $\circles$ are all congruent the bound improves to
\begin{equation} \label{eqnew-uni}
\tilde O\left(\frac{n^{3/2}}{k^{11/6}} + \frac{n^{3/2}q^{1/2}}{k^3} +
\frac{n}{k}   \right).
\end{equation}
In particular, if $q=O(1)$, the number of such points is
\[ \tilde O\left(\frac{n^{3/2}}{k^{11/6}} + \frac{n}{k}   \right). \]
\end{corollary}
\begin{theProof}{\!\!}
Let $m$ be the number of points incident to at least $k$ circles of $\circles$, and observe that these points determine
at least $mk$ incidences with the circles of $\circles$.  Comparing this
lower bound with the upper bound in Theorem \ref{th:arb} (for (a)),
or in Theorem \ref{th:uni} (for (b)), the claims follow.
\end{theProof}

\noindent{\bf Remarks.}
(1) It is interesting to compare the bounds in (\ref{eqnew}) and (\ref{eqnew-uni})
with the various recent bounds on incidences between points and lines in three
dimensions~\cite{EKS11,GK10,GK11}. In all of them the threshold value
$m=\Theta(n^{3/2})$ plays a significant role. Specifically:
(i) The number of joints in a set of $n$ lines in $\RR^3$ is
$O(n^{3/2})$, a bound tight in the worst case \cite{GK10}.
(ii) If no plane contains more than $\sqrt{n}$ lines, the number of points
incident to at least $k\ge 3$ lines is $O(n^{3/2}/k^2)$ \cite{GK11}.
(iii) A related bound where $m=n^{3/2}$ is a threshold value, under
different assumptions, is given in \cite{EKS11}.
The bounds in (\ref{eqnew}) and (\ref{eqnew-uni}) are somewhat weaker
(because of the extra small factors hidden in the $\tilde O(\cdot)$ notation,
the rather restrictive constraints on $q$, and the constraint $k\ge k_0$)
but they belong to the same class of results.
It would be interesting to understand how general this phenomenon is; for
example, does it also show up in incidences with other classes of curves in $\RR^3$?
We tend to conjecture that this is the case, under reasonable assumptions
concerning those curves.  Similar threshold phenomena should exist in higher
dimensions. ``Extrapolating" from the results of \cite{KSS10,Quil10}, these
thresholds should be at $m=n^{d/(d-1)}$.
\medskip

\noindent
(2) The bounds can be slightly tightened by using Theorem~\ref{th:noEps} or
Theorem~\ref{th:UnitNoEps} instead of Theorem~\ref{th:arb} or Theorem~\ref{th:uni},
respectively, but we leave these slight improvements to the interested reader.

\paragraph{Similar triangles.}
Another application of Theorem \ref{th:arb} (or rather of Theorem \ref{th:noEps})
is an improved bound on the
number of triangles spanned by a set $\pts$ of $t$ points in $\RR^3$ and
similar to a given triangle $\Delta$.  Let $F(\pts,\Delta)$ be the number of triangles spanned by $\pts$ that are similar to $\Delta$, and let $F(t)$ be the maximum of $F(\pts,\Delta)$ as $\pts$ ranges over all sets of $t$ points and $\Delta$ ranges over all triangles. We then have:

\begin{theorem}
$$
F(t)=O(t^{15/7})=O(t^{2.143}).
$$
\end{theorem}

\begin{theProof}{\!\!}
Let $\pts$ be a set of $t$ points in $\RR^3$ and let $\Delta=uvw$ be a
given triangle. Suppose that $pqr$ is a similar copy of $\Delta$, where
$p,q,r\in \pts$. If $p$ corresponds to $u$ and $q$ to $v,$ then $r$ has
to lie on a circle $c_{pq}$ that is orthogonal to the segment $pq$, whose
center lies at a fixed point on this segment, and whose radius is proportional
to $|pq|$. Thus, the number of possible candidates for the point $r$, for
$p,q$ fixed, is exactly the number of incidences between $\pts$ and $c_{pq}$.
There are $2\binom{t}{2} = t(t-1)$ such circles, and no circle arises more
than twice in this manner.  It follows that $F(t)$ is bounded by twice the
number of incidences between the $t$ points of $\pts$ and the $t(t-1)$
circles $c_{pq}$. We now apply Theorem \ref{th:noEps} with $m=t$ and
$n=t(t-1)$. (The theorem applies for these values, which satisfy
$m\approx n^{1/2}$, much smaller than the threshold $n^{3/2}$; in fact,
$m$ lies in the second range $[n^{1/3},n^{4/5}]$.)
It remains to show that the expression \eqref{mainBound} is $O(t^{15/7})$.

The first term of \eqref{mainBound} is $O(t^{15/7})$. To control the
remaining terms, it suffices to show that at most
$O\Big(\big(\frac{n^3}{m^2} \big)^{3/7} \Big)=O(t^{12/7})$ of the circles
lie on a common plane or sphere.  In fact, we claim that at most $O(t)$
circles can lie on a common plane or sphere. Indeed, let $\Pi$ be a plane.
Then for any circle $c_{pq}$ contained in $\Pi$, $pq$ must be orthogonal
to $\Pi,$ pass through the center of $c_{pq}$, and each of $p$ and $q$
must lie at a fixed distance from $\Pi$ (the distances are determined
by the triangle $\Delta$ and by the radius of $c_{pq}$).
This implies that each point of $\pts$ can generate at most two
circles on $\Pi$. The argument for cosphericality is essentially the
same. The only difference is that one point of $\pts$ may lie at the
center of the given sphere $\sigma$, and then it can determine up to
$2(t-1)$ distinct circles on $\sigma$. Still, the number of circles on
$\sigma$ is $O(t)$.
As noted above, this completes the proof of the theorem.
\end{theProof}

As already mentioned in the introduction, this slightly improves a previous
bound of $O^*(n^{58/27})$ in \cite{AS11} (see also \cite{aaps07}).

\paragraph{Acknowledgements.}
We would like to thank Ji\v{r}\'i Matou\v{s}ek for organizing the
\emph{Micro-Workshop on Algebraic Methods in Discrete Geometry}
in Z\"{u}rich in May 2012, where the authors made significant progress on the above proof, and we gratefully acknowledge support by the ERC grant
DISCONV that facilitated the workshop.  We would also like to thank Noam
Solomon for helpful discussions about algebraic geometry.  The first
two authors would like to thank the people that participated in an algebraic
geometry reading seminar with them, and for several useful discussions.
They are Roel Apfelbaum, Dan Halperin, Haim Kaplan, Manjish Pal, Orit Raz,
and Shakhar Smorodinsky.  The third author would like to thank J\'ozsef
Solymosi for introducing this problem to him.

We would also like to thank Ji\v{r}\'i Matou\v{s}ek and Zuzana Safernov\'a
for pointing out errors in a previous draft of this manuscript.

\end{document}